\documentclass[11pt,reqno]{amsart}
\usepackage{amsmath, amsthm, amssymb}
\usepackage{enumitem}
\usepackage[foot]{amsaddr}

\usepackage{fullpage}
\usepackage[many]{tcolorbox}
\usepackage{xcolor}
\usepackage{wasysym}
\usepackage{mathtools}
\usepackage{scalerel}
\usepackage{bbm}
\usepackage[all]{xy}
\usepackage{mathrsfs}
\usepackage{adjustbox}
\usepackage{tensor}

\usepackage{tikz}
\usepackage{adjustbox}
\usetikzlibrary{arrows,backgrounds,patterns.meta}
\usetikzlibrary{positioning,shadings,cd}
\usetikzlibrary{shapes}
\usetikzlibrary{backgrounds}
\usetikzlibrary{decorations,decorations.pathreplacing,decorations.markings,decorations.pathmorphing}
\tikzstyle{snake}=[decorate, decoration={snake, segment length=1mm, amplitude=.5mm}]
\usetikzlibrary{fit,calc,through}
\usetikzlibrary{external}

\tikzset{super thick/.style={line width=3pt}}
\tikzstyle{far>}=[decoration={markings, mark=at position 0.75 with {\arrow{>}}}, postaction={decorate}]
\tikzstyle{mid>}=[decoration={markings, mark=at position 0.55 with {\arrow{>}}}, postaction={decorate}]
\tikzstyle{mid<}=[decoration={markings, mark=at position 0.55 with {\arrow{<}}}, postaction={decorate}]
\tikzset{super thick/.style={line width=3pt}}
\tikzstyle{far>}=[decoration={markings, mark=at position 0.75 with {\arrow{>}}}, postaction={decorate}]
\tikzstyle{mid>}=[decoration={markings, mark=at position 0.55 with {\arrow{>}}}, postaction={decorate}]
\tikzstyle{mid<}=[decoration={markings, mark=at position 0.55 with {\arrow{<}}}, postaction={decorate}]
\tikzstyle{knot}=[preaction={super thick, white, draw}]
\tikzstyle{coupon}=[draw, very thick, rectangle, rounded corners=5pt]
\tikzset{Rightarrow/.style={double equal sign distance,>={Implies},->},
triplecd/.style={-,preaction={draw,Rightarrow}},
quadruplecd/.style={preaction={draw,Rightarrow,
shorten >=0pt
},
shorten >=1pt,
-,double,double
distance=0.2pt}}
\tikzset{
    tripleline/.style args={[#1] in [#2] in [#3]}{
        #1,preaction={preaction={draw,#3},draw,#2}
    }
}
\tikzstyle{triple}=[tripleline={[line width=.15mm,black] in
      [line width=.7mm,white] in
      [line width=1mm,black]}] 
\tikzset{
    quadrupleline/.style args={[#1] in [#2] in [#3] in [#4]}{
        #1,preaction={preaction={preaction={draw,#4},draw,#3}, draw,#2}
    }
}
\tikzstyle{quadruple}=[quadrupleline={[line width=.3mm,white] in
      [line width=.6mm,black] in
      [line width=1.2mm,white] in
      [line width=1.5mm,black]}]

\definecolor{violet}{RGB}{148,0,211}
\definecolor{DarkGreen}{RGB}{0,150,0}
\definecolor{rufous}{HTML}{A81C07}

\usepackage[plainpages=false,hypertexnames=false,pdfpagelabels]{hyperref}
\definecolor{medium-blue}{rgb}{0,0,.8}
\hypersetup{colorlinks, linkcolor={purple}, citecolor={medium-blue}, urlcolor={medium-blue}}
\newcommand{\arxiv}[1]{\href{http://arxiv.org/abs/#1}{\tt arXiv:\nolinkurl{#1}}}
\newcommand{\arXiv}[1]{\href{http://arxiv.org/abs/#1}{\tt arXiv:\nolinkurl{#1}}}

\newcommand{\googlebooks}[1]{(preview at \href{https://books.google.com/books?id=#1}{google books})}

\DeclareMathOperator{\End}{End}

\DeclareMathOperator{\ev}{ev}

\DeclareMathOperator{\id}{id}

\DeclareMathOperator{\ind}{ind}

\DeclareMathOperator{\Irr}{Irr}

\DeclareMathOperator{\op}{op}

\DeclareMathOperator{\Vect}{Vec}
\DeclareMathOperator{\Hilb}{Hilb}

\DeclareMathOperator{\Mod}{Mod}
\DeclareMathOperator{\Alg}{Alg}
\DeclareMathOperator{\CMod}{C^*Mod}
\DeclareMathOperator{\WMod}{W^*Mod}
\DeclareMathOperator{\HilbMod}{HilbMod}
\DeclareMathOperator{\CAlg}{C^*\!Alg}
\DeclareMathOperator{\WAlg}{W^*\!Alg}

\newcommand{\bb}[1]{\mathbf{#1}}
\newcommand{\bbover}[1]{\mathbf{\overline{#1}}}

\def\semicolon{;}
\def\applytolist#1{
    \expandafter\def\csname multi#1\endcsname##1{
        \def\multiack{##1}\ifx\multiack\semicolon
            \def\next{\relax}
        \else
            \csname #1\endcsname{##1}
            \def\next{\csname multi#1\endcsname}
        \fi
        \next}
    \csname multi#1\endcsname}

\def\calc#1{\expandafter\def\csname c#1\endcsname{{\mathcal #1}}}
\applytolist{calc}QWERTYUIOPLKJHGFDSAZXCVBNM;
\def\bbc#1{\expandafter\def\csname bb#1\endcsname{{\mathbb #1}}}
\applytolist{bbc}QWERTYUIOPLKJHGFDSAZXCVBNM;
\def\bfc#1{\expandafter\def\csname bf#1\endcsname{{\mathbf #1}}}
\applytolist{bfc}QWERTYUIOPLKJHGFDSAZXCVBNM;
\def\sfc#1{\expandafter\def\csname s#1\endcsname{{\sf #1}}}
\applytolist{sfc}QWERTYUIOPLKJHGFDSAZXCVBNM;
\def\fc#1{\expandafter\def\csname f#1\endcsname{{\mathfrak #1}}}
\applytolist{fc}QWERTYUIOPLKJHGFDSAZXCVBNM;
\def\rmc#1{\expandafter\def\csname rm#1\endcsname{{\mathrm #1}}}
\applytolist{rmc}QWERTYUIOPLKJHGFDSAZXCVBNM;
\def\scrc#1{\expandafter\def\csname scr#1\endcsname{{\mathscr #1}}}
\applytolist{scrc}QWERTYUIOPLKJHGFDSAZXCVBNM;

\numberwithin{equation}{section}

\theoremstyle{plain}
\newtheorem{thm}{Theorem}[section]
\newtheorem{defn}[thm]{Definition}

\newtheorem{cor}[thm]{Corollary}
\newtheorem{prop}[thm]{Proposition}
\newtheorem{ex}[thm]{Example}

\newtheorem{lemma}[thm]{Lemma}

\newtheorem{rmk}[thm]{Remark}

\newtheorem{const}[thm]{Construction}
\newtheorem*{thm*}{Theorem}
\newtheorem*{rmk*}{Remark}

\newtheorem*{claim*}{Claim}

\newtheorem{thmalpha}{Theorem}

\title{Ostrik's Theorem for Operator Algebras in a Unitary Tensor Category}
\author[E.~Sherman]{Emeline Sherman}

\begin{document}

\begin{abstract}
We define Hilbert modules for $\rmC^*$ and $\rmW^*$ algebra objects in the ind-completion of a unitary tensor category.
We then prove an operator algebraic Ostrik’s theorem in an ind-category: the 2-category of $\rmC^*$ or $\rmW^*$ algebras objects, their modules, and adjointable intertwiners is $\dagger$-equivalent to the 2-categories of $\rmC^*$ or $\rmW^*$ module categories.
\end{abstract}

\maketitle

\section{Introduction}
Ostrik's theorem \cite{Ost03} shows that for a fusion category $\cC$, every $\cC$ module category is a category of modules over an algebra in $\cC$.
This theorem is best stated as an equivalence of 2-categories, $\Alg_{\cC}$ and $\Mod_{\cC}$ respectively:
\[
  \left\{ \begin{array}{c}
    \text{algebra objects in }\cC \\
    \text{bimodules in }\cC \\
    \text{intertwiners}
  \end{array}\right\}
  \text{ and }
  \left\{ \begin{array}{c}
    \text{cyclic } \cC \text{ module categories} \\
    \text{module functors} \\
    \text{natural transformations}
  \end{array}\right\}.
\]
From an algebra object $A$ in $\cC$, its category of modules $\Mod_{\cC}(A)$, the Cauchy completion\footnote{For linear categories, Cauchy completion is equivalent to being closed under direct sums and sub-objects.} of the free modules, is a module category.
On the other hand, from a $\cC$ module category $\cM$, the internal hom functor produces an algebra object $\underline \End_{\cC}(m)$ for some $m \in \cM$.
This recovers the Morita equivalence class of an algebra object.
When $\cM$ is cyclic, i.e. has a basepoint $n$, $\underline \End_{\cC}(n)$ is a particular algebra object in a Morita equivalence class.
For example, $\underline \End_{\cC}(A) \cong A$.
Ostrik notes in \cite[Remark~9]{Ost03} that is is quite easy to extend this theorem to the case of $\cC$ a semisimple rigid tensor category with simple unit, using ind-objects in the ind-completion of $\cC$.

For $\cC$ a unitary tensor category, i.e., a semisimple rigid $\rmC^*$-tensor category with simple unit,
C. Jones and Penneys developed a well-behaved theory of operator algebra objects in \cite{JP17}.
The completeness conditions for a $*$-algebra object in $\ind$-$\cC$ are certified using Ostrik's theorem.
More precisely, a $*$-algebra object $\bfA$ in $\ind$-$\cC$ is a $\rmC^*$ or $\rmW^*$ algebra object if and only if its corresponding cyclic $\cC$ module category $\Mod_{\cC}(\bfA)$ is $\rmC^*$ or $\rmW^*$, respectively.
This provides an equivalence of 1-categories, where the morphisms of $\rmC^*$ algebra objects are algebra natural transformations, and the morphisms of cyclic $\cC$ module $\rmC^*$ categories are cyclic module functors, both of course with the appropriate properties.
In this paper, we develop the necessary technology to extend this to an equivalence of $\rmC^*/\rmW^*$ 2-categories.

Recall from \cite{Lan95,BLM04,MT05,Bla05} that a Hilbert $\rmC^*$-module for an honest $\rmC^*$-algebra $A$ is a vector space with an $A$-valued inner product, under which it is complete.
The choices of 1-morphisms as Hilbert $\rmC^*$-bimodules (often called right correspondences in the literature) and of 2-morphisms between Hilbert $\rmC^*$-bimodules as adjointable intertwiners assemble a 2-category $\CAlg$.
Similarly, the 2-category $\WAlg$ consists of $\rm W^*$ algebra objects, their bimodules, and normal intertwiners.
We develop these module objects in $\ind$-$\cC$ for a unitary tensor category $\cC$.
This motivates our main result:
\begin{thmalpha} \label{thm:A}
There exist $\dagger$-equivalences of 2-categories $\CAlg_{\cC} \cong \CMod_{\cC}$ and $\WAlg_{\cC} \cong \WMod_{\cC}$.
In more detail,
\[
  \left\{ \begin{array}{c}
    \rmC^* \text{ algebra objects} \\
    \text{Hilbert } \rmC^* \text{ bimodule objects} \\
    \text{adjointable intertwiners}
  \end{array}\right\}
  \cong
  \left\{ \begin{array}{c}
    \text{cyclic } \cC \text{ module } \rmC^* \text{ categories} \\
    \text{dagger module functors} \\
    \text{uniformly bounded natural transformations}
  \end{array}\right\}
  \]
  and
  \[
  \left\{ \begin{array}{c}
    \rmW^* \text{ algebra objects} \\
    \text{Hilbert } \rmW^* \text{ bimodule objects} \\
    \text{normal intertwiners}
  \end{array}\right\}
  \cong
  \left\{ \begin{array}{c}
    \text{cyclic } \cC \text{ module } \rmW^* \text{ categories} \\
    \text{normal dagger module functors} \\
    \text{uniformly bounded natural transformations}
  \end{array}\right\}.
\]
\end{thmalpha}
Note that we assume all $\rmW^*$-algebra objects are tracial.\footnote{We are currently developing Tomita-Takesaki theory for operator algebra objects in a unitary tensor category.
}

\begin{rmk*}
Theorem \ref{thm:A} is a generalization of \cite[Theorem~7.4]{HNP24}.
More precisely, our result is in an arbitrary unitary tensor category $\cC$ what their result is in $\Hilb_{\rm fd}$.
See Remark \ref{rmk:HNPgeneralized} for further discussion.
\end{rmk*}
Theorem \ref{thm:A} is also similar to the well-known 2-categorical equivalence of standard separable $Q$-systems and algebra objects in $\cC$ with appropriate morphisms; see a sketch in \cite[Example~3.39]{CHPJP22}.
See also \cite[Theorem~3.21]{Ver22}, proven in full detail in \cite[Appendix~6.2]{Ver22}.
Whereas the known result deals with (operator) algebra objects in $\cC$, our result addresses the ind completion of $\cC$.
We will elaborate on this connection in future work.

In order to define $\WAlg_{\cC}$ properly, we need to obtain a result similar to the identification of the two forms of bimodules for von Neumann algebras.
A Hilbert $\rmW^*$ module for a $\rmW^*$ algebra is a Hilbert $\rmC^*$ module that is compatible with the predual.
On the other hand, we can describe a Hilbert space module for a von Neumann algebra as a normal representation on Hilbert space.
Just as $\rmW^*$ algebra and von Neumann algebras are equivalent \cite{Sak98}, it is well-known that Hilbert $\rmW^*$ modules and Hilbert space modules form $\dagger$-equivalent categories \cite[8.5.39]{BLM04}, \cite[1.4]{DH94}, \cite[Theorem~2.2]{BDH88}.
From the first, take the relative product with the standard form of the appropriate von Neumann algebra.
From the second, take the left bounded vectors (when working with right modules).
With this in mind, we prove the following result about Hilbert $\rmW^*$ module objects and Hilbert space module objects.
\begin{thmalpha} \label{thm:B}
The functors $$- \otimes_{\bb M} L^2 \bb M \colon \WMod_{\bb M} \to \HilbMod_{\bb M}$$ and $$D(-) \coloneq \HilbMod_{\bb M}(L^2 \bb M, -) \colon \HilbMod_{\bb M} \to \WMod_{\bb M}$$ induce a $\dagger$-equivalence of $\rmW^*$ categories.
\end{thmalpha}

\subsection*{Acknowledgements}
We are very grateful to David Penneys for pointing us towards this problem and his guidance in the development of the theory.
In particular, Definition \ref{defn:HilbCModObj} is due to him, and much of the theory of Section \ref{section:CModObjs} was developed jointly.
We also thank Gregory Faurot, Giovanni Ferrer, and Brett Hungar for helpful conversations.
The author was supported by the Department of Mathematics at Ohio State University, a fellowship from the Quantum Graduate Interdisciplinary Program funded by OSU's NSF NRT grant 2244045, and Penneys' NSF DMS grant 2554723.

\subsection*{AI use statement}
No LLMs were used in the preparation of this manuscript, nor for any of the results it contains.
\section{Background}
In this section, we review the relevant material of \cite{JP17}.
Let $\cC$ be a unitary tensor category, i.e. a semisimple rigid $\rmC^*$ tensor category with simple unit.
We work in the ind-completion and unitary ind-completion (see \cite{NY16}) of $\cC$.
Our model for the first is $\Vect(\cC)$, the category of linear functors $\cC^{\op} \to \Vect$ with linear natural transformations.
For the second, we use $\Hilb(\cC)$, the category of linear dagger functors $\cC^{\op} \to \Hilb$ and uniformly bounded natural transformations, i.e. $\alpha$ such that $\|\alpha\| \coloneq \sup \limits_{c \in \cC} \|\alpha_c\| < \infty$, where $\|\alpha_c\|$ is the operator norm.
The way to think about objects in these categories is as $\cC$-graded vector or Hilbert spaces.
The interested reader can take $\cC$ to be $\Hilb_{\rm fd}$, upon which $\Vect(\cC) \cong \Vect$ and $\Hilb(\cC) \cong \Hilb$.
Doing so recovers the usual definitions of operator algebras and their modules.

These categories inherit a monoidal structure from $\cC$ via Day convolution \cite{Day70}.
For $\bb X, \bb Y \in \Vect(\cC)$, $$\bb X \otimes \bb Y(c) \coloneq \bigoplus \limits_{a, b \in \Irr(\cC)} \bb X(a) \otimes \cC(c, a \otimes b) \otimes \bb Y(b)$$
There is a similar construction for $\Hilb(\cC)$.
For each component Hilbert space $\bb H \otimes \bb K(c)$, the inner product is scaled by $d_a^{-1} d_b^{-1}$, where $d_c$ is the quantum dimension of an object $c \in \cC$.
We write the $c$ component as $\bigoplus \limits_{a, b \in \Irr(\cC)}\bb H(a) \otimes_a \cC(c, a \otimes b) \otimes_b \bb K(b)$.
Conveniently, elements of the form $\eta \otimes \alpha \otimes \xi$ span a dense subset in each component.
Using this, we describe the $c$ component of the associator $\alpha^{\Vect(\cC)}_{\bb X, \bb Y, \bb Z} \colon \bb X \otimes (\bb Y \otimes \bb Z) \to (\bb X \otimes \bb Y) \otimes \bb Z$ by
$$\eta \otimes \alpha \otimes (\xi \otimes \beta \otimes \zeta) \mapsto \sum U^{\gamma, \delta}_{\alpha, \beta} (\eta \otimes \gamma \otimes \xi) \otimes \delta \otimes \zeta$$
\begin{center}
\begin{tikzpicture}[baseline={(current bounding box.center)}, thick]
   \filldraw [black] (0,0) circle (0pt) coordinate (a);
   \filldraw [black] ([xshift=0.6cm, yshift=0.5cm]a) circle (0pt) coordinate (b);
   \draw (a) to [out=0, in=-90] node[pos=0.5, below] {$\bb b$} (b);
   \draw ([xshift=-1.25cm, yshift=1.75cm]a) node[draw=black, fill=white, very thick, rounded corners, baseline=center] (c){$\eta$};
   \draw ([xshift=1.25cm]c) node[draw=black, fill=white, very thick, rounded corners, baseline=center] (d){$\xi$};
   \draw ([xshift=1.25cm]d) node[draw=black, fill=white, very thick, rounded corners, baseline=center] (e){$\zeta$};
   \node[above] at (a) {$\alpha$};
   \node[above] at (b) {$\beta$};
   \draw (a) to [out=180, in=-90] node[pos=0.5, left] {$\bb a$} (c.south);
   \draw (b) to [out=180, in=-90] node[pos=0.5, left] {$\bb d$} (d.south);
   \draw (b) to [out=0, in=-90] node[pos=0.5, right] {$\bb e$} (e.south);
   \draw (a) --node[pos=1, below] {$\bb c$}++ (0, -0.5);
   \draw (c.north) --node[pos=1, above] {$\bb X$}++ (0, 0.5) coordinate (f);
   \draw (d.north) to node[pos=1, above] {$\bb Y$} ([xshift=1.25cm]f);
   \draw (e.north) to node[pos=1, above] {$\bb Z$} ([xshift=2.5cm]f);

   \node[] (f) at ([xshift=1.5cm, yshift=-0.4cm]e) {$\mapsto \sum U_{\alpha, \beta}^{\gamma, \delta}$};

   \filldraw [black] ([xshift=2.75cm, yshift=-1.35cm]f) circle (0pt) coordinate (a);
   \filldraw [black] ([xshift=-0.6cm, yshift=0.5cm]a) circle (0pt) coordinate (b);
   \draw (a) to [out=180, in=-90] node[pos=0.5, below] {$\bb f$} (b);
   \draw ([xshift=-1.25cm, yshift=1.75cm]a) node[draw=black, fill=white, very thick, rounded corners, baseline=center] (c){$\eta$};
   \draw ([xshift=1.25cm]c) node[draw=black, fill=white, very thick, rounded corners, baseline=center] (d){$\xi$};
   \draw ([xshift=1.25cm]d) node[draw=black, fill=white, very thick, rounded corners, baseline=center] (e){$\zeta$};
   \node[above] at (a) {$\delta$};
   \node[above] at (b) {$\gamma$};
   \draw (a) to [out=0, in=-90] node[pos=0.5, right] {$\bb e$} (e.south);
   \draw (b) to [out=0, in=-90] node[pos=0.5, right] {$\bb d$} (d.south);
   \draw (b) to [out=180, in=-90] node[pos=0.5, left] {$\bb a$} (c.south);
   \draw (a) --node[pos=1, below] {$\bb c$}++ (0, -0.5);
   \draw (c.north) --node[pos=1, above] {$\bb X$}++ (0, 0.5) coordinate (f);
   \draw (d.north) to node[pos=1, above] {$\bb Y$} ([xshift=1.25cm]f);
   \draw (e.north) to node[pos=1, above] {$\bb Z$} ([xshift=2.5cm]f);
\end{tikzpicture}
\end{center}
where the sum is taken over $f \in \Irr(\cC)$, an ONB of $\gamma$'s for $\cC(f, a \otimes d)$, and an ONB of $\delta$'s for $\cC(c, f \otimes e)$.
Again, one can similarly write down the associator for $\Hilb(\cC)$.
We clarify a few details of the above diagram.
The scalars $U^{\gamma, \delta}_{\alpha, \beta}$ are entries in matrices for $F$-symbols in $\cC$.
The objects of the form $\bb c$ are functors $\cC(-, c)$, where by $\bbover c$ we mean $\cC(-, c^\vee)$.\footnote{Objects in $\Vect(\cC)$ and $\Hilb(\cC)$ will in general be written in bold, and elements in $\cC$ will not.}\textsuperscript{,}\footnote{Since objects in $\Vect(\cC)$ and $\Hilb(\cC)$ are too ``large" to be dualizable, passing from $\cC$ to $\Vect(\cC)$ and $\Hilb(\cC)$ forgets the unitary dual functor $(\cdot)^\vee$.}
We call these objects \emph{compact}.
These objects inherit \emph{all} of the structure and property of $\cC$, since they are the image of the Yoneda embedding of $\cC$.
Now using Yoneda's lemma, we (uniquely) identify $\eta \in \bb X(a)$, for example, with a natural transformation $\eta\colon \bb a \to \bb X$.
The unitors follow from a basic semisimple category computation.
We suppress the associator and unitor natural isomorphisms whenever possible.

$\Vect(\cC)$ and $\Hilb(\cC)$ also admit an involutive structure, the conjugate.
On a category $\cD$, an involutive structure is an anti-linear functor $\bar{\cdot} \colon \cD \to \cD$ with a natural isomorphisms $\varphi \colon \id_{\cD} \to \bar{\bar{\cdot}}$, $\nu_{a, b} \colon \bar{a} \otimes \bar{b} \to \overline{b \otimes a}$ for all $a, b \in \cD$, and an isomorphism $r \colon 1_{\cD} \to \overline{1_{\cD}}$.
These satisfy some relations, importantly
\begin{equation} \label{eq:assoc}
\nu_{a, c \otimes b} \circ (\id_{\overline{a}} \otimes \nu_{b, c}) = \nu_{b \otimes a, c} \circ (\nu_{a, b} \otimes \id_{\overline c}) 
\end{equation}
\begin{equation} \label{eq:decomp}
\varphi_{a \otimes b} = \overline{\nu_{b,a}} \circ \nu_{\bar{a}, \bar{b}} \circ (\varphi_a \otimes \varphi_b)
\end{equation}
For $\bb X \in \Vect(\cC)$, $\bbover X (c) \coloneq \bbover X(c^\vee)$.
In the interest of clarity, $\bbover X$ is a linear functor, but $\overline{\cdot}$ is an \emph{anti-linear} functor.
We direct the interested reader to \cite{JP17} for more details, in particular for the definition of $\nu$ in $\Vect(\cC)$.
As with unitors and associators, we will freely suppress involutive data.

\subsection{Operator Algebra Objects in Unitary Tensor Categories}
\begin{defn}
An algebra object in $\Vect(\cC)$ is a triple $\bfA, m_{\bfA}, \iota_{\bfA}$, the latter two of which we represent as diagrams of the form
\[
\begin{tikzpicture}[baseline={(current bounding box.center)}, thick]
   \filldraw [black] (0,0) circle (0pt) coordinate (a);
   \draw ([xshift=-0.5cm, yshift=-0.5cm]a) to [out=90, in=180] node[pos=0, below] {$\bb A$} (a);
   \draw ([xshift=0.5cm, yshift=-0.5cm]a) to [out=90, in=0] node[pos=0, below] {$\bb A$} (a);
   \draw (a) --node[pos=1, above] {$\bb A$}++ (0,0.5);
\end{tikzpicture}
\quad \text{and} \quad\quad
\begin{tikzpicture}[baseline={(current bounding box.center)}, thick]
   \filldraw [black] (0,0) circle (2pt) coordinate (a);
   \draw (a) --node[pos=1, above] {$\bb A$}++ (0,0.5);
\end{tikzpicture}
\quad \text{respectively}\]
satisfying
\begin{center}
\begin{tikzpicture}[baseline={(current bounding box.center)}, thick]
   \filldraw [black] (0,0) circle (0pt) coordinate (a);
   \filldraw [black] ([xshift=0.5cm, yshift=0.5cm]a) circle (0pt) coordinate (b);
   \draw ([xshift=-0.5cm, yshift=-0.5cm]a) to [out=90, in=180] node[pos=0, below] {$\bb A$} (a);
   \draw ([xshift=0.5cm, yshift=-0.5cm]a) to [out=90, in=0] node[pos=0, below] {$\bb A$} (a);
   \draw ([xshift=1.25cm, yshift=-0.5cm]a) to [out=90, in=0] node[pos=0, below] {$\bb A$} (b);
   \draw (a) to [out=90, in=180] (b);
   \draw (b) --node[pos=1, above] {$\bb A$} ++ (0, 0.5);
   
   \node[] (c) at ([xshift=1.25cm, yshift=-0.25cm]b) {$=$};
   
   \filldraw [black] ([xshift=1.75cm, yshift=-0.25cm]c) circle (0pt) coordinate (a);
   \filldraw [black] ([xshift=-0.5cm, yshift=0.5cm]a) circle (0pt) coordinate (b);
   \draw ([xshift=-0.5cm, yshift=-0.5cm]a) to [out=90, in=180] node[pos=0, below] {$\bb A$} (a);
   \draw ([xshift=0.5cm, yshift=-0.5cm]a) to [out=90, in=0] node[pos=0, below] {$\bb A$} (a);
   \draw ([xshift=-1.25cm, yshift=-0.5cm]a) to [out=90, in=180] node[pos=0, below] {$\bb A$} (b);
   \draw (a) to [out=90, in=0] (b);
   \draw (b) --node[pos=1, above] {$\bb A$} ++ (0, 0.5);
   
   \node[] (c) at ([xshift=2.25cm, yshift=-0.25cm]b) {and};
   
   \filldraw [black] ([xshift=1.75cm, yshift=-0.25cm]c) circle (2pt) coordinate (a);
   \draw ([xshift=0.5cm, yshift=-0.5cm]a) --node[pos=0, below] {$\bb A$} node[pos=1, above] {$\bb A$}++ (0,1.5);
   \filldraw [black] ([xshift=0.5cm, yshift=0.5cm]a) circle (0pt) coordinate (b);
   \draw (a) to [out=90, in=180] (b);
   
   \node[] (c) at ([xshift=0.5cm, yshift=-0.25cm]b) {$=$};
   
   \filldraw [black] ([xshift=1cm, yshift=-0.25cm]c) circle (0pt) coordinate (a);
   \draw ([xshift=-0.5cm, yshift=-0.5cm]a) --node[pos=0, below] {$\bb A$} node[pos=1, above] {$\bb A$}++ (0,1.5);
   
   \node[] (c) at ([xshift=0cm, yshift=0.25cm]a) {$=$};
   
   \filldraw [black] ([xshift=1cm, yshift=-0.25cm]c) circle (2pt) coordinate (a);
   \draw ([xshift=-0.5cm, yshift=-0.5cm]a) --node[pos=0, below] {$\bb A$} node[pos=1, above] {$\bb A$}++ (0,1.5);
   \filldraw [black] ([xshift=-0.5cm, yshift=0.5cm]a) circle (0pt) coordinate (b);
   \draw (a) to [out=90, in=0] (b);
\end{tikzpicture}
\end{center}
A $*$-algebra object $\bb A$ is an algebra object equipped with a $*$-structure, i.e., a conjugate linear natural transformation $j \colon \bb A \to \bb A$ such that $j^2 = \id_{\bb A}$ and $j$ is compatible with the monoidal and involutive structures of $\bb A$ (as an object in $\Vect(\cC)$).
Equivalently, we can describe this $*$-structure as an involutive algebra natural isomorphism $\sigma \colon \bb A \to \bbover A$ such that $\overline \sigma \circ \sigma = \varphi$.
\end{defn}

Next, we define module objects for algebra objects.

\begin{defn}
For an algebra object $\bb A \in \Vect(\cC)$, a right $\bb A$ module object is $\bb X \in \Vect(\cC)$ together with a natural transformation $\rho_{\bb A} \colon \bb X \otimes \bb A \to \bb X$ such that

\begin{center}
\begin{tikzpicture}[baseline={(current bounding box.center)},thick]
\draw[red,very thick] (0,0)--node[pos=0,black,below]{$\bb X$}++(0, 2);
\draw (1, 0) to [bend left=-45] node[pos=0,black,below]{$\bb A$} (0, 1);
\draw (1.4, 0) to [bend left=-45] node[pos=0,black,below]{$\bb A$}(0, 1.4);
\node[] at (2,1) {$=$};
\draw[red,very thick] (2.8,0)--node[pos=0,black,below]{$\bb X$}++(0, 2);
\draw (3.8, 0.5) to [bend left=-45] node[pos=0.7,black,above right]{$\bb A$} (2.8, 1.4);
\draw (3.4, 0) to [bend left=45] node[pos=0,black,below]{$\bb A$} (3.8, 0.5);
\draw (4.2, 0) to [bend right=45] node[pos=0,black,below]{$\bb A$} (3.8, 0.5);
\filldraw [black] (3.8, 0.5) circle (0pt);

\node[] at (5.5, 1) {and};

\draw [red, very thick] (7,0) --node[pos=0,black,below]{$\bb X$}++ (0,2);
\draw (8,0.5) to [bend left = -45] node[pos=0.7,black,above right]{$\bb A$} (7, 1.4);
\filldraw [black] (8, 0.5) circle (2pt);
\node[] at (9,1) {$=$};
\draw[red, very thick] (9.8,0) -- node[pos=0,black,below]{$\bb X$}++(0,2);
\end{tikzpicture}
\end{center}

We can similarly define a left $\bb A$ module object.
Here, we denote the action $\bb A \otimes \bb X \to \bb X$ by $\lambda_{\bb A}$.
When $\bb X$ is a left or right $\bb A$ module object, we write $_{\bb A} \bb X$ or $\bb X_{\bb A}$, respectively.
A natural transformation $\alpha \colon \bb X_{\bb A} \to \bb Y_{\bb A}$ is right $\bb A$-linear if the following diagram holds.

\begin{center}
\begin{tikzpicture}[baseline={(current bounding box.center)}, thick]
\draw (0,0) node[draw=black,fill=white,very thick,rounded corners,baseline=center] (a){$\alpha$};
\draw [red, very thick] (a.south) -- node[pos=1, black, below]{$\bb X$} coordinate (b) ++(0, -0.75);
\draw [red, very thick] (a.north) -- node[pos=1, black, above]{$\bb Y$} ++(0, 0.75);
\draw ([xshift=0.5cm, yshift=-0.75cm]a.south) to [out=90, in=0]  node[pos=0, below]{$\bb A$} (b);

\node[] (c) at ([xshift=1cm]a) {$=$};

\draw ([xshift=1cm]c) node[draw=black,fill=white,very thick,rounded corners,baseline=center] (d){$\alpha$};
\draw [red, very thick] (d.south) -- node[pos=1, black, below]{$\bb X$} ++(0, -0.75);
\draw [red, very thick] (d.north) -- node[pos=1, black, above]{$\bb Y$} coordinate (e) ++(0, 0.75);
\draw ([xshift=0.75cm, yshift=-0.75cm]d.south) to [out=90, in=0]  node[pos=0, below]{$\bb A$} (e);
\end{tikzpicture}
\end{center}
\end{defn}

We will soon see that right $\bb A$-module objects and right $\bfA$ module maps form a left $\cC$ module category, denoted $\cM_{\bb A}$.

\begin{defn} A (left) $\cC$ module category $\cM$ consists of the following data.
	\begin{itemize}
	\item A bilinear bifunctor $-\otimes_{\cM}- \colon \cC \otimes \cM \to \cM$
	\item Unitor and associator natural isomorphisms for $\otimes_{\cM}$, i.e. for all $m \in \cM$, $a, b, c, d \in \cC$
		\begin{itemize}
		\item $\lambda \colon 1_{\cC} \otimes m \to m$
		\item $\rho \colon m \otimes 1_{\cC} \to m$
		\item $\alpha_{c,d,m} \colon c \otimes (d \otimes m) \to (c \otimes_{\cC} d) \otimes m$
		\end{itemize}
	\end{itemize}
satisfying
\begin{align*}
\alpha_{a, b \otimes_{\cC} c, m} \circ \alpha_{a, b, c \otimes m} &= (\alpha_{a, b, c} \otimes 1_m) \circ \alpha_{a \otimes_{\cC} b, c, m} \circ (1_a \otimes \alpha_{b, c, m})\\
(\rho_c \otimes 1_m) \circ \alpha_{1_{\cC}, c, m} = 1_c \otimes \lambda_m &\text{ and }
(\lambda_c \otimes 1_m) \circ \alpha_{1_{\cC},c,m} = \lambda_{c \otimes m}
\end{align*}

We say $\cM$ is cyclic, written $\cM_m$, if $\cM$ comes equipped with a choice of basepoint $m \in \cM$.
A $\cC$ module category $\cM$ is a $\cC$ module $\dagger$-category if $\cM$ is a dagger category, the dagger structures on $\cC$ and $\cM$ respect $\otimes_{\cM}$ (i.e. $(\phi \otimes \psi)^* = \phi^* \otimes \psi^*$), and $\alpha$, $\lambda$, and $\rho$ are unitary isomorphisms.
\end{defn}

To define operator algebra objects, \cite{JP17} used the unitary 1-categorical version of Ostrik's theorem.
Here, the morphisms of algebra objects are algebra natural transformations, and the morphisms of cyclic $\cC$ module categories are cyclic module functors, defined below.

\begin{defn} A module functor $\cM \to \cN$ for $\cC$ module categories is a functor together with natural isomorphisms $c \otimes F(m) \stackrel{\omega_{c, m}}{\to} F(c \otimes m)$ satisfying associativity and unitality diagrams
$$\xymatrix{c \otimes (d \otimes F(m)) \ar[r]^{1_c \otimes \omega_{d, m}} \ar[d]_\alpha & c \otimes F(d \otimes m) \ar[r]^{\omega_{c, d \otimes m}} & F(c \otimes (d \otimes m)) \ar[d]^{F(\alpha)} \\ (c \otimes d) \otimes F(m) \ar[rr]^{\omega_{c \otimes d, m}} & & F((c \otimes d) \otimes m) }$$

$$\xymatrix{1_{\cC} \otimes F(m) \ar[rr]^{\omega_{c, m}} \ar[dr]_{\lambda_{F(m)}} && F(1_{\cC} \otimes m) \ar[dl]^{F(\lambda_m)} \\ & F(m)}$$

When $\cM_m$ and $\cN_n$ are cyclic, $F$ is a \emph{cyclic} module functor if there exists an isomorphism $F(m) \cong n$.
If $\cM$ and $\cN$ are $\cC$ module dagger categories, we require $F$ to be a dagger functor.
In this case, for a cyclic module (dagger) functor the isomorphism $F(m) \cong n$ must be unitary.
\end{defn}

Consider the 2-category of cyclic $\cC$ module categories, cyclic module functors, and natural transformations.
This category truncates to a 1-category.
Specifically, \cite[Lemma~3.6]{HPT23b} tells us that all 2-morphisms in this 2-category are invertible.
Thus we can consider the 1-morphisms to be equivalence classes under natural transformations of cyclic module functors.
The equivalence of this 1-category with the category of algebra objects in $\Vect(\cC)$ and algebra natural transformations is \cite[Remark~9]{Ost03}.
In the case of both algebra and $*$-algebra objects, we sketch the equivalence on objects for the convenience of the reader.
To do so, we observe that by Yoneda's lemma, the categories of ($*$-)algebra objects with ($*$-preserving) algebra natural transformations and (dagger) lax monoidal functors with (dagger) monoidal natural transformations are equivalent.

For an algebra object $\bb A$, consider the cyclic $\cC$ module category $\cM_{\bb A}$ with basepoint $\bb A$.
When $\bb A$ is a $*$-algebra object, we equip $\cM_{\bb A}$ with a dagger structure by leveraging $$\cM_{\bb A}(c \otimes \bb A, d \otimes \bb A) \cong \Vect(\cC)(\bb c, \bb d \otimes \bb A) \cong \bb A(\overline d \otimes c)$$
The first isomorphism is \cite[Lemma~3.2]{Ost03}, and the second holds by Frobenius reciprocity.
For $f \in \bb A(\overline b \otimes a)$, $$f^* \coloneq (\bb A(\id_a \otimes \varphi_b) \circ \bb A(\nu_{a, \overline b} \circ j_{\overline b \otimes a})(f) \in \bb A(\overline a \otimes b)$$
On the other hand, let $\cM_m$ be a cyclic $\cC$ module category.
We define a lax monoidal functor $\bb A$ with components $\bb A(c) \coloneq \cM_m(c \otimes m, m)$.
The tensorator $\mu_{a, b}$ for $\bb A$ is given by $f \otimes g \mapsto f \circ (\id_a \otimes g)$.
The unitor is given by $\id_m \circ \lambda \in \cM(1_{\cC} \otimes m, m)$.
When $\cM_m$ is a dagger category, $\bb A$ is a $*$-algebra object via the diagram:
\begin{center}
\begin{tikzpicture}[baseline={(current bounding box.center)},thick]
\draw (0, 0) node[draw=black,fill=white,very thick,rounded corners,baseline=center] (a){$f^*$};
\draw ([xshift=-0.2cm]a.north) to [out=90, in=0] ([xshift=-0.4cm, yshift=0.2cm]a.north) coordinate (d);
\draw (d) to [out=180, in=90] ([xshift=-0.6cm]a.north) coordinate (b);
\draw ([xshift=0.2cm]a.north) --node[pos=1, above] {$m$}++ (0, 0.5);
\draw (a.south) --node[pos=1, below] {$m$}++ (0, -0.5) coordinate (c);
\draw (b) to node[pos=1, below] {$a$} ([xshift=-0.6cm]c);

\node[] at ([xshift=-2cm]a) {$j_a(f) \coloneq$};
\node[] at ([xshift=2cm]a) {for $f \in \bb A(a)$};
\end{tikzpicture}
\end{center}

Since $\cC$ (hence $\Vect(\cC)$ and $\Hilb(\cC)$) have direct sums, so too does $\cM$.
Thus, $\cM$ is a cyclic
\footnote{By remembering the generating object, one can harmlessly pass between the category of free modules for $\bfA$, with object $\bf c \otimes \bf A$ for $c \in \cC$, and the category of right modules for $\bfA$, the Cauchy completion of the free $\bfA$ modules.
For example, a cyclic $\cC$ module category is a $\rmC^*$ or $\rmW^*$ category \emph{if and only if} its Cauchy completion is.
We will often work with the category of free modules for convenience.}
$\cC$ module $\rmC^*$ or $\rmW^*$ category if for all $c \in \cC$, $\cM(c \otimes m, c \otimes m)$ is a $\rmC^*$ or von Neumann algebra, respectively.
From this, we can define operator algebra objects in a unitary tensor category.

\begin{defn}
For an algebra object $\bb A$ in $\Vect(\cC)$, $\bb A$ is a $\rmC^*$ or $\rmW^*$ algebra object if $\bb A(\overline c \otimes c)$ is a $\rmC^*$ or von Neumann algebra for all $c \in \cC$.
\end{defn}

We provide some examples and non-examples of $\rmW^*$ (hence $\rmC^*$) algebra objects.

\begin{ex}
Let $c \in \cC$.
As $\cC$ is a $\rmW^*$ category, $\bbover c \otimes \bb c$ is a $\rmW^*$ algebra object.
\end{ex}

\begin{ex}
Consider $\bb H \in \Hilb(\cC)$.
We take the cyclic $\cC$ module category it generates $\cM_{\bb H}$.
As $\Hilb(\cC)$ is a $\rmW^*$ category, so too is $\cM_{\bb H}$.
So the algebra object corresponding to this module category, $B(\bb H)$, is a $\rmW^*$ algebra object.

For a non-example, consider $\bb X \in \Vect(\cC)$.
We denote the algebra object corresponding to the module category $\cM_{\bb X}$ by $ \underline{\End}_{\cC}(\bb X)$.
This will \emph{not} be an operator algebra object in general.
\end{ex}

Note that for $\bb X \in \Vect(\cC)$, a right $\bb A$ module object $\bb X$ can be interpreted as an algebra natural transformation $\bb A^{\op} \to \underline \End_{\cC}(\bb X)$.
Here, $\bb A^{\op}$ has multiplication given by
$$\xymatrix{\bb A^{\op} (a) \otimes \bb A^{\op} (b) = \bb A(a^\vee) \otimes \bb A(b^\vee) \ar[r]^<<<<\beta & \bb A(b^\vee) \otimes \bb A(a^\vee) \ar[r]^<<<<\mu & \bb A(b^\vee \otimes a^\vee) \cong \bb A^{\op} (a \otimes_{\cC^{\op}} b})$$
Similarly, for $\bb H \in \Hilb(\cC)$, we can view $\bb H$ as a right $\bb A$ module object via an ($*$-)algebra natural transformation $\bb A^{\op} \to B(\bb H)$.
We also call this a ($*$-)representation of $\bb A$.

The last piece of theory from \cite{JP17} that we will need is a suitable notion of $L^2 \bb M$ for a tracial $\rmW^*$ algebra object $\bb M$.
A $\rmC^*$ algebra object $\bb A$ comes with a left and a right $\bb A(1_{\cC})$-valued inner product on each component:
\[{}_a\langle f , g \rangle \coloneq
\begin{tikzpicture}[baseline={(current bounding box.center)},thick]
\draw (0, 0) node[draw=black,fill=white,very thick,rounded corners,baseline=center] (a){$f$};
\draw ([xshift=1cm]a) node[draw=black,fill=white,very thick,rounded corners,baseline=center] (b){$j_a(g)$};
\filldraw ([xshift=0.5cm, yshift=0.75cm]a) circle (0pt) coordinate (c);
\filldraw ([xshift=0.5cm, yshift=-0.75cm]a) circle (0pt) coordinate (d);
\draw (a.north) to [out=90, in=180] node[pos=0.4, left] {$\bb A$} (c);
\draw (b.north) to [out=90, in=0] node[pos=0.4, right] {$\bb A$} (c);
\draw (a.south) to [out=-90, in=180] node[pos=0.4, left] {$\bbover a$} (d);
\draw (b.south) to [out=-90, in=0] node[pos=0.4, right] {$\bb a$} (d);
\draw (c) --node[pos=1, above] {$\bb A$}++ (0, 0.4);
\end{tikzpicture}
\quad and \quad
\langle f | g \rangle_a \coloneq
\begin{tikzpicture}[baseline={(current bounding box.center)},thick]
\draw (0, 0) node[draw=black,fill=white,very thick,rounded corners,baseline=center] (a){$j_a(f)$};
\draw ([xshift=1cm]a) node[draw=black,fill=white,very thick,rounded corners,baseline=center] (b){$g$};
\filldraw ([xshift=0.5cm, yshift=0.75cm]a) circle (0pt) coordinate (c);
\filldraw ([xshift=0.5cm, yshift=-0.75cm]a) circle (0pt) coordinate (d);
\draw (a.north) to [out=90, in=180] node[pos=0.4, left] {$\bb A$} (c);
\draw (b.north) to [out=90, in=0] node[pos=0.4, right] {$\bb A$} (c);
\draw (a.south) to [out=-90, in=180] node[pos=0.4, left] {$\bb a$} (d);
\draw (b.south) to [out=-90, in=0] node[pos=0.4, right] {$\bbover a$} (d);
\draw (c) --node[pos=1, above] {$\bb A$}++ (0, 0.4);
\end{tikzpicture}
\]
If $\bb A$ admits a (faithful) state $\varphi$, by which we mean $\bb A(1_{\cC})$ admits a state, this turns the $\bb A(1_{\cC})$-valued inner product into a genuine inner product by post-composition.
Completing each component $\bb A(c)$ with respect to the first inner product gives ${}_\varphi L^2\bb A \in \Hilb(\cC)$, and with respect to the second gives $L^2\bb A_\varphi \in \Hilb(\cC)$.
Note that the left ${}_\varphi L^2 \bb A$ has a right $\bb A$ action, and the right $L^2 \bb A_\varphi$ has a left $\bb A$ action.
These are given by right or left multiplication, respectively, on $\widehat{\bb A}$, which is the component-wise image of $\bb A$ under the GNS construction.
Observe that $\widehat{\bb A}$ is component-wise dense in $L^2\bb A$.
We can alternatively view $L^2\bb A_\varphi$ as having $c$ component $\bb A(c) \otimes_{\bb A(1_{\cC})} L^2(\bb A(1_{\cC}), \varphi)$, and similarly for ${}_\varphi L^2 \bb A$.
For clarity, this has inner product $\langle \eta_1 \otimes \xi_1 | \eta_2 \otimes \xi_2 \rangle \coloneq \langle \xi_1 | \langle \varphi(\eta_1 | \eta_2 \rangle_a) \xi_2 \rangle_{L^2(\bb A(1_{\cC}), \varphi)}$

The left and right versions of $L^2 \bb A$ are conjugate isometric.
Thus, when we take the state to be a trace $\tau$, the left and right constructions agree.
So we have left \emph{and} right $\bb A$ actions on $L^2 \bb A$.
If $\bb M$ is a $\rmW^*$ algebra object with (faithful normal) trace $\tau$, these actions are normal and commute (up to an associator in $\Hilb(\cC)$).
By normal, we mean their actions component-wise are normal.
We also call this a normal $*$-representation of $\bb M$.

\section{Module Objects for $\rmC^*$ Algebra Objects}\label{section:CModObjs}
The first task we undertake is defining $\CAlg_{\cC}$, the 2-category of $\rmC^*$ algebra objects, Hilbert $\rmC^*$ bimodule objects, and adjointable intertwiners.

\begin{defn}\label{defn:HilbCModObj}
For a $\rmC^*$ algebra object $\bb A \in \Vect(\cC)$ and a right $\bb A$ module object $\bb X$, we say that $\bb X$ when taken together with a form $\langle \cdot |\cdot \rangle_\bb A \colon \bbover X \otimes \bb X \to \bb A$ is a right Hilbert $\rmC^*$ module object when the following conditions are satisfied.

\begin{enumerate}[label=\textup{(H\arabic*)}]
\item \label{H1} Right linearity$\colon$
\begin{tikzpicture}[baseline={(current bounding box.center)},thick]
\draw (0, 0) node[draw=black,fill=white,very thick,rounded corners,baseline=center] (a){$\langle\cdot|\cdot\rangle_{\mathbf{A}}$};
\draw [red, very thick] ([xshift=-0.3cm]a.south) --node[pos=1, below, black] {$\bbover X$}++ (0, -0.75) coordinate (b);
\draw [red, very thick] ([xshift=0.3cm]a.south) to node[pos=1, below, black] {$\bb X$} coordinate (c) ([xshift=0.6cm]b);
\draw ([xshift=1.1cm]b) to [out=90, in=0] node[pos=0, below] {$\bb A$} (c);
\draw (a.north) --node[pos=1, above] {$\bb A$}++ (0, 0.5);

\node[] (d) at ([xshift=0.5cm]a.east) {$=$};

\draw ([xshift=1cm, yshift=-0.25cm]d) node[draw=black,fill=white,very thick,rounded corners,baseline=center] (a){$\langle\cdot|\cdot\rangle_{\mathbf{A}}$};
\draw [red, very thick] ([xshift=-0.3cm]a.south) --node[pos=1, below, black] {$\bbover X$}++ (0, -0.5) coordinate (b);
\draw [red, very thick] ([xshift=0.3cm]a.south) to node[pos=1, below, black] {$\bb X$} ([xshift=0.6cm]b);
\draw (a.north) --node[pos=1, above] {$\bb A$}  coordinate (c)++ (0, 0.75);
\draw ([xshift=1.5cm]b) to [out=90, in=0] node[pos=0, below] {$\bb A$} (c);
\end{tikzpicture}

\item \label{H2} Anti-symmetry$\colon$
\begin{tikzpicture}[baseline={(current bounding box.center)},thick]
\draw (0, 0) node[draw=black,fill=white,very thick,rounded corners,baseline=center] (a){$\langle\cdot|\cdot\rangle_{\mathbf{A}}$};
\draw [red, very thick] ([xshift=-0.3cm]a.south) to node[pos=1, black, below]{$\bbover X$}++(0, -0.5) coordinate (b);
\draw [red, very thick] ([xshift=0.3cm]a.south) to node[pos=1, black, below]{$\bbover X$} ([xshift=0.6cm]b);
\draw (a.north) --node[pos=1, black, above]{$\bb A$}++(0, 0.5);

\node[] at ([xshift=1cm]a) {$=$};

\draw ([xshift=2cm]a) node[draw=black,fill=white,very thick,rounded corners,baseline=center] (b){$\overline{\langle\cdot|\cdot\rangle_{\mathbf{A}}}$};
\draw [red, very thick] ([xshift=-0.05cm,yshift=-0.75cm]b.south) --++(0, 0.75);
\draw [red, very thick] ([xshift=0.05cm,yshift=-0.75cm]b.south) --node[pos=0.45, black, right]{$\overline{\bbover X \otimes \bb X}$}++(0, 0.75);
\draw (b.north) --node[pos=0.5, black, right] {$\bbover A$}++(0,0.5);

\draw ([yshift=1.25cm]b) node[draw=black,fill=white,very thick,rounded corners,baseline=center] (c){$\sigma^{-1}_{\bb A}$};
\draw (c.north)--node[pos=1, black, above] {$\bb A$}++(0, 0.5);

\draw ([yshift=-1.45cm]b) node[draw=black,fill=white,very thick,rounded corners,baseline=center] (d){$\nu_{\bbover X, \bb X}$};
\draw [red, very thick] ([xshift=-0.3cm, yshift=-1cm]d.south) coordinate (z) --node[pos=0, black, below] {$\bbover X$}++(0, 1);
\draw [red,very thick] ([xshift=0.3cm, yshift=-0.25cm]d.south) --++(0, 0.25);
\node[] at ([xshift=0.75cm, yshift=0cm]d.south) {$\overline{\bbover X}$};
\draw [red, very thick] ([xshift=0.3cm]d.south) to node[pos=1, black, below] {$\bb X$} ([xshift=0.6cm]z);
\draw ([xshift=0.3cm, yshift=-0.5cm]d.south) node[draw=black,fill=white,very thick,rounded corners,baseline=center] (e){$\varphi_{\bb X}$};
\end{tikzpicture}
\item \label{H3} Positive-definiteness: for all $c \in \cC$ and $\eta \in \bb X(c)$, the following is a positive element of $\bb A(\overline{c} \otimes c)$, and 0 iff $\eta = 0$.

\begin{center}
\begin{tikzpicture}[baseline={(current bounding box.center)},thick]
\draw (0, 0) node[draw=black,fill=white,very thick,rounded corners,baseline=center] (a){$\langle \cdot | \cdot \rangle_{\bb A}$};
\draw [red, very thick] ([xshift=-0.3cm, yshift=-0.75cm]a.south) --node[pos=0.5, left, black] {$\bbover X$}++(0, 0.75);
\draw [red, very thick] ([xshift=0.3cm, yshift=-0.75cm]a.south) --node[pos=0.5, right, black] {$\bb X$}++(0, 0.75);
\draw (a.north) --node[pos=1, above, black] {$\bb A$}++(0, 0.5);

\draw ([xshift=-0.3cm, yshift=-1.35cm]a) node[draw=black,fill=white,very thick,rounded corners,baseline=center] (b){$\overline \eta$};
\draw (b.south) --node[pos=1, below]{$\bbover c$}++(0, -0.5) coordinate (d);

\draw ([xshift=0.3cm, yshift=-1.35cm]a) node[draw=black,fill=white,very thick,rounded corners,baseline=center] (c){$\eta$};
\draw (c.south) to node[pos=1, below]{$\bb c$} ([xshift=0.6cm]d);
\end{tikzpicture}
\end{center}
We denote this diagram by $\langle \xi | \xi \rangle_{\bb A}$.

\item \label{H4} Completeness$\colon$ $\bb X(c)$ is a right $\bb A(1_\cC)$ Hilbert $\rmC^*$ module for all $c \in \cC$, i.e. $\bb X(c)$ is complete with respect to the norm $\|\eta\|_{\bb X(c)} \coloneq \|\langle \eta | \eta \rangle_{\bb A}^\frac{1}{2}\|_{\bb A(1_{\cC})}$.
\end{enumerate}
When the form $\langle \cdot | \cdot \rangle$ satisfies $\ref{H1}-\ref{H4}$, we call it a right $\bb A$-inner product.
We can similarly define a left Hilbert $\rmC^*$ module object.
Here, we write the left $\bb A$-valued inner product as $_{\bb A} \langle \cdot, \cdot \rangle$.
\end{defn}

If desired, one can swap \ref{H1} for the following condition (see Proposition \ref{prop:RightLeftSwap}).
\begin{enumerate}[label=\textup{(H\arabic*')}]
\item \label{H1'} Left anti-linearity$\colon$
\begin{tikzpicture}[baseline={(current bounding box.center)},thick]
\draw (0, 0) node[draw=black,fill=white,very thick,rounded corners,baseline=center] (a){$\langle\cdot|\cdot\rangle_{\mathbf{A}}$};
\draw [red, very thick] ([xshift=-0.3cm]a.south) --node[pos=1, below, black] {$\bbover X$} coordinate (c) ++ (0, -1) coordinate (b);
\draw [red, very thick] ([xshift=0.3cm]a.south) to node[pos=1, below, black] {$\bb X$} ([xshift=0.6cm]b);
\draw ([xshift=-0.5cm]b) to [out=90, in=180] node[pos=0, below] {$\bbover A$} (c);
\draw (a.north) --node[pos=1, above] {$\bb A$}++ (0, 1);

\node[] (d) at ([xshift=0.5cm]a.east) {$=$};

\draw ([xshift=2cm, yshift=-0.25cm]d) node[draw=black,fill=white,very thick,rounded corners,baseline=center] (a){$\langle\cdot|\cdot\rangle_{\mathbf{A}}$};
\draw [red, very thick] ([xshift=-0.3cm]a.south) --node[pos=1, below, black] {$\bbover X$}++ (0, -0.75) coordinate (b);
\draw [red, very thick] ([xshift=0.3cm]a.south) to node[pos=1, below, black] {$\bb X$} ([xshift=0.6cm]b);
\draw (a.north) --node[pos=1, above] {$\bb A$}  coordinate (c)++ (0, 1.25);
\draw ([xshift=-1.25cm, yshift=-0.25cm]a) node[draw=black,fill=white,very thick,rounded corners,baseline=center] (e){$\sigma^{-1}_A$};
\draw (e.north) to [out=90, in=180] node[pos=0.7, above left] {$\bb A$} (c);
\draw (e.south) to node[pos=1, below] {$\bbover A$} ([xshift=-0.95cm]b);
\end{tikzpicture}

\[\text{where the vertex } \begin{tikzpicture}[baseline={(current bounding box.center)}, thick]
\draw[red,very thick] (0,0)coordinate(x)--node[pos=0,black,below]{$\bbover X$}coordinate(md)++(0, 1);
\draw([xshift=-0.5cm]x)to[out=90,in=180]node[pos=0,black,below]{$\bbover A$}(md);
\end{tikzpicture}
\text{ represents the morphism }
\xymatrix{\rho_{\bbover A} \coloneq \bbover A \otimes \bbover X \ar[r]^>>>>{\nu_{\bb A, \bb X}} & \overline{\bb X \otimes \bb A} \ar[r]^>>>>{\overline{\lambda_{\bb A}}} & \bbover X \,.}\]
\end{enumerate}

There are two useful weakenings of Definition \ref{defn:HilbCModObj}.
We note these now, for later use in Constructions \ref{const:completion} and \ref{const:relprod}.
\begin{defn}
We say $\bb X$ is a right $\bb A$ pre-Hilbert $\rmC^*$ module object if it satisfies \ref{H1}, \ref{H2}, \ref{H3}, but not \ref{H4}.
Similarly, if a right $\bb A$-valued inner product $\langle \cdot | \cdot \rangle_{\bb A}$ satisfies \ref{H1}, \ref{H2}, but not \ref{H3}, we call $\langle \cdot | \cdot \rangle_{\bb A}$ a pre-inner product if it instead satisfies
\begin{enumerate}
\item [(H3')] \label{H3'} Positivity$\colon$ for all $c \in \cC$ and $\xi \in \bb X(c)$, $\langle \xi | \xi \rangle_{\bb A} \geq 0$.
\end{enumerate}
\end{defn}

Observe the following.
\begin{prop} \label{prop:RightLeftSwap}
\ref{H1} and \ref{H2} imply \ref{H1'}.
Similarly, \ref{H1'} and \ref{H2} imply \ref{H1}
\end{prop}

\begin{proof}
We show the first statement holds.
Equalities not otherwise justified follow by naturality or associativity (\ref{eq:assoc}) of $\nu$, or naturality of the algebra multiplication.

\begin{center}
\begin{tikzpicture}[baseline={(current bounding box.center)},thick]
\draw (0, 0) node[draw=black,fill=white,very thick,rounded corners,baseline=center] (a){$\langle \cdot | \cdot \rangle_{\bb A}$};
\draw (a.north) --node[pos=1, above] {$\bb A$}++(0, 1);
\draw [red, very thick] ([xshift=-0.3cm]a.south) --node[black, pos=1, below] {$\bbover X$} coordinate(b) ++(0, -1);
\draw [red, very thick] ([xshift=0.3cm]a.south) --node[black, pos=1, below] {$\bb X$}++(0, -1);
\draw ([xshift=-0.8cm, yshift=-1cm]a.south) to [out=90, in=180] node[pos=0, below] {$\bb A$} (b);

\node[] (c) at ([xshift=1.25cm]a){$\stackrel{\ref{H2}}{=}$};

\draw ([xshift=1.25cm]c) node[draw=black,fill=white,very thick,rounded corners,baseline=center,inner xsep=0.5cm] (d){$\nu$};
\draw [red, very thick] ([xshift=-0.3cm]d.south) --node[black, left, pos=0.5] {$\bbover X$} ++(0, -0.75);
\draw [red, very thick] ([xshift=0.3cm]d.south) --node[black, right, pos=0.5] {$\overline{\bbover X}$} ++(0, -0.75);
\draw [red, very thick] ([xshift=-0.05cm]d.north) --++(0, 0.75);
\draw [red, very thick] ([xshift=0.05cm]d.north) --node[black, right, pos=0.5] {$\overline{\bbover X \otimes \bb X}$} ++(0, 0.75);
\draw ([yshift=1.25cm]d) node[draw=black,fill=white,very thick,rounded corners,baseline=center] (e){$\langle \cdot | \cdot \rangle_{\bb A}$};
\draw(e.north) --node[right, pos=0.5] {$\bbover A$} ++(0, 0.75);
\draw ([yshift=1.25cm]e) node[draw=black,fill=white,very thick,rounded corners,baseline=center] (f){$\sigma^{-1}$};
\draw (f.north) --node[above, pos=1] {$\bb A$} ++(0, 0.5);
\draw ([xshift=-0.4cm,yshift=-1.25cm]d) node[draw=black,fill=white,very thick,rounded corners,baseline=center] (g){$\overline{\lambda_{\bb A}}$};
\draw ([xshift=0.4cm,yshift=-1.25cm]d) node[draw=black,fill=white,very thick,rounded corners,baseline=center] (h){$\varphi$};
\draw [red, very thick]([xshift=-0.05cm]g.south) --node[black, left, pos=0.5] {$\bbover{X \otimes A}$}++(0, -0.75);
\draw ([xshift=0.05cm]g.south) --++(0, -0.75);
\draw ([yshift=-1.25cm]g) node[draw=black,fill=white,very thick,rounded corners,baseline=center, inner xsep =0.3cm] (i){$\nu$};
\draw ([xshift=-0.2cm]i.south) --node[pos=1, below] {$\bbover A$} ++(0, -0.5) coordinate (j);
\draw [red, very thick] ([xshift=0.2cm]i.south) --node[black, pos=1, below] {$\bbover X$} ++(0, -0.5);
\draw [red, very thick] (h.south) to node[black, pos=1, below] {$\bb X$} ([xshift=1cm]j);

\node[] (s) at ([xshift=1.75cm]d){$\stackrel{\ref{H1}}{=}$};

\draw ([xshift=2cm]s) node[draw=black,fill=white,very thick,rounded corners,baseline=center,inner xsep=0.5cm] (t){$\overline{\langle \cdot | \cdot \rangle_{\bb A} \otimes 1_{\bb A}}$};
\draw ([xshift=0.1cm]t.south) --node[black, right, pos=0.45] {$\overline{\bbover X \otimes \bb X \otimes \bb A}$} ++(0, -0.75);
\draw [red, very thick] ([xshift=-0.1cm]t.south) --++(0, -0.75);
\draw [red, very thick] (t.south) --++(0, -0.75);
\draw ([xshift=0.05cm]t.north) --node[black, right, pos=0.4] {$\overline{\bb A \otimes \bb A}$} ++(0, 0.75);
\draw ([xshift=-0.05cm]t.north) --++(0, 0.75);
\draw ([yshift=1.25cm]t) node[draw=black,fill=white,very thick,rounded corners,baseline=center] (u){$\overline{m_{\bb A}}$};
\draw(u.north) --node[right, pos=0.5] {$\bbover A$} ++(0, 0.75);
\draw ([yshift=1.25cm]u) node[draw=black,fill=white,very thick,rounded corners,baseline=center] (v){$\sigma^{-1}$};
\draw (v.north) --node[above, pos=1] {$\bb A$} ++(0, 0.5);
\draw ([yshift=-1.25cm]t) node[draw=black,fill=white,very thick,rounded corners,baseline=center, inner xsep=0.5cm] (w){$\nu$};
\draw [red, very thick]([xshift=-0.45cm]w.south) --node[black, left, pos=0.5] {$\bbover{X \otimes A}$}++(0, -0.75);
\draw ([xshift=-0.35cm]w.south) --++(0, -0.75);
\draw [red, very thick] ([xshift=0.4cm]w.south) --node[pos=0.45, black, right] {$\overline{\bbover X}$} ++(0, -0.75);
\draw ([xshift=0.4cm,yshift=-1.25cm]w) node[draw=black,fill=white,very thick,rounded corners,baseline=center] (x){$\varphi$};
\draw ([xshift=-0.4cm,yshift=-1.2cm]w) node[draw=black,fill=white,very thick,rounded corners,baseline=center, inner xsep=0.3cm] (y){$\nu$};
\draw [red, very thick] ([xshift=0.2cm]y.south) --node[black, pos=1, below] {$\bbover X$} ++(0, -0.5) coordinate (z);
\draw ([xshift=-0.2cm]y.south) --node[pos=1, below] {$\bbover A$} ++(0, -0.5);
\draw [red, very thick] ([xshift=0.2cm]i.south) --node[black, pos=1, below] {$\bbover X$} ++(0, -0.5);
\draw [red, very thick] (x.south) to node[black, pos=1, below] {$\bb X$} ([xshift=0.6cm]z);

\node[] at ([xshift=2.25cm]t){$=$};

\draw ([xshift=3.5cm]t) node[draw=black,fill=white,very thick,rounded corners,baseline=center, inner xsep=0.5cm] (h){$\nu$};
\draw ([xshift=0.05cm]h.north) --node[black, right, pos=0.4] {$\overline{\bb A \otimes \bb A}$} ++(0, 0.75);
\draw ([xshift=-0.05cm]h.north) --++(0, 0.75);
\draw ([yshift=1.25cm]h) node[draw=black,fill=white,very thick,rounded corners,baseline=center] (i){$\overline{m_{\bb A}}$};
\draw(i.north) --node[right, pos=0.5] {$\bbover A$} ++(0, 0.75);
\draw ([yshift=1.25cm]i) node[draw=black,fill=white,very thick,rounded corners,baseline=center] (j){$\sigma^{-1}$};
\draw (j.north) --node[above, pos=1] {$\bb A$} ++(0, 0.5);
\draw ([xshift=-0.4cm]h.south) --node[pos=1, below] {$\bbover A$} ++(0, -3.25) coordinate (l);
\draw ([xshift=0.4cm]h.south) --node[pos=0.4, right] {$\bbover A$} ++(0, -0.75);
\draw ([xshift=0.4cm, yshift=-1.25cm]h) node[draw=black,fill=white,very thick,rounded corners,baseline=center] (j){$\overline{\langle \cdot | \cdot \rangle_{\bb A}}$};
\draw [red, very thick] ([xshift=-0.05cm]j.south) --++ (0, -0.75);
\draw [red, very thick] ([xshift=0.05cm]j.south) --node[pos=0.4, right, black] {$\overline{\bbover X \otimes \bb X}$} ++ (0, -0.75);
\draw ([yshift=-1.25cm]j) node[draw=black,fill=white,very thick,rounded corners,baseline=center, inner xsep = 0.4cm] (k){$\nu$};
\draw [red, very thick] ([xshift=-0.2cm]k.south) to node[black, pos=1, below] {$\bbover X$} ([xshift=0.6cm]l);
\draw [red, very thick] ([xshift=0.2cm]k.south) to node[pos=1, below, black] {$\bb X$} ([xshift=1.0cm]l);

\node[] (m) at ([xshift=1.75cm]h){$=$};

\filldraw ([xshift=1.75cm, yshift=1.25cm]m) circle (0pt) coordinate (t);
\draw ([xshift=1cm]m) node[draw=black,fill=white,very thick,rounded corners,baseline=center] (u){$\sigma^{-1}_{\bb A}$};
\draw ([xshift=2.5cm]m) node[draw=black,fill=white,very thick,rounded corners,baseline=center] (v){$\sigma^{-1}_{\bb A}$};
\draw (u.north) to [out=90, in=180] node[pos=0.3, above left] {$\bb A$} (t);
\draw (v.north) to [out=90, in=0] node[pos=0.3, above right] {$\bb A$} (t);
\draw (t) --node[pos=1, above] {$\bb A$} ++(0, 0.75);
\draw (u.south) --node[pos=1, below] {$\bbover A$} ++(0, -3.25) coordinate (y);
\draw (v.south) --node[pos=0.3, right] {$\bbover A$} ++(0, -0.75);
\draw ([yshift=-1.25cm]v) node[draw=black,fill=white,very thick,rounded corners,baseline=center] (w){$\overline{\langle \cdot | \cdot \rangle_{\bb A}}$};
\draw [red, very thick] ([xshift=-0.05cm]w.south) --++ (0, -0.75);
\draw [red, very thick] ([xshift=0.05cm]w.south) --node[pos=0.4, right, black] {$\overline{\bbover X \otimes \bb X}$} ++ (0, -0.75);
\draw ([yshift=-1.25cm]w) node[draw=black,fill=white,very thick,rounded corners,baseline=center, inner xsep = 0.4cm] (x){$\nu$};
\draw [red, very thick] ([xshift=-0.2cm]x.south) to node[black, pos=1, below] {$\bbover X$} ([xshift=1.3cm]y);
\draw [red, very thick] ([xshift=0.2cm]x.south) to node[pos=1, below, black] {$\bb X$} ([xshift=1.7cm]y);
\end{tikzpicture}
\end{center}
An application of $\ref{H2}$ completes the proof.
\end{proof}

We supply some examples of Hilbert $\rmC^*$ module objects.
Our first example is the expected one.
\begin{ex} \label{ex:overself}
The $\rmC^*$ algebra object $\bb A$ is a right (and left) $\bb A$ Hilbert $\rmC^*$ module object.
We define the right $\bb A$-valued inner product as follows:
\[
\begin{tikzpicture}[baseline={(current bounding box.center)}, thick]
\draw (0, 0) node[draw=black,fill=white,very thick,rounded corners,baseline=center] (a){$\sigma_A^{-1}$};
\draw (a.south) --node[pos=1, below]{$\bbover A$} ++(0, -0.25);
\draw(a.north) to [out=90,in=180] node[pos=0.5, above left]{$\bb A$} ([xshift=0.5cm, yshift=0.5cm]a.north);
\filldraw ([xshift=0.5cm, yshift=0.5cm]a.north) circle (0pt);
\draw ([xshift=1cm, yshift=-0.25cm]a.south) to [out=90, in=0] node[pos=0, below]{$\bb A$} ([xshift=0.5cm, yshift=0.5cm]a.north);
\draw([xshift=0.5cm, yshift=0.5cm]a.north) -- node[pos=1, above]{$\bb A$} ++(0, 0.5);
\end{tikzpicture}
\] This is right-linear by associativity of $\bb A$.
The $\bb A$-valued inner product is anti-symmetric:

\begin{center}
\begin{tikzpicture}[baseline={(current bounding box.center)}, thick]
\draw (0, 0) node[draw=black,fill=white,very thick,rounded corners,baseline=center] (a){$\overline{\sigma^{-1}_{\bb A} \otimes id_{\bb A}}$};
\draw ([xshift=-0.05cm,yshift=-0.75cm]a.south) --++(0, 0.75);
\draw ([xshift=0.05cm,yshift=-0.75cm]a.south) --node[pos=0.45, black, right]{$\overline{\bbover A \otimes \bb A}$}++(0, 0.75);
\draw ([xshift=-0.05cm]a.north) --++(0,0.75);
\draw ([xshift=0.05cm]a.north) --node[pos=0.5, black, right] {$\overline{\bbover A \otimes \bb A}$}++(0,0.75);

\draw ([yshift=1.5cm]a) node[draw=black,fill=white,very thick,rounded corners,baseline=center] (b){$\overline{m_{\bb A}}$};
\draw (b.north) --node[pos=0.4, right]{$\bbover A$}++(0, 0.75);

\draw ([yshift=2.75cm]a) node[draw=black,fill=white,very thick,rounded corners,baseline=center] (b){$\sigma^{-1}_{\bb A}$};
\draw (b.north)--node[pos=1, black, above] {$\bb A$}++(0, 0.25);

\draw ([yshift=-1.45cm]a) node[draw=black,fill=white,very thick,rounded corners,baseline=center] (c){$\nu_{\bbover A, \bb A}$};
\draw ([xshift=-0.3cm, yshift=-1cm]c.south) --node[pos=0, black, below] {$\bbover A$}++(0, 1);
\draw ([xshift=0.3cm, yshift=-0.25cm]c.south) --++(0, 0.25);
\node[] at ([xshift=0.75cm, yshift=0cm]c.south) {$\overline{\bbover A}$};

\draw ([xshift=0.3cm, yshift=-0.5cm]c.south) node[draw=black,fill=white,very thick,rounded corners,baseline=center] (d){$\varphi_{\bb A}$};
\draw ([yshift=-0.25cm]d.south) --node[pos=0, black, below] {$\bb A$}++(0, 0.25);

\node[] (e) at ([xshift=1.75cm]a){$=$};

\draw ([xshift=1.75cm]e) node[draw=black,fill=white,very thick,rounded corners,baseline=center] (j){$\overline{\sigma^{-1}_{\bb A} \otimes id_{\bb A}}$};
\draw ([xshift=-0.05cm]j.north) -- ++(0, 0.75);
\draw ([xshift=0.05cm]j.north) --node[pos=0.5, right]{$\overline{\bb A \otimes \bb A}$} ++(0, 0.75);

\draw ([yshift=1.5cm]j) node[draw=black,fill=white,very thick,rounded corners,baseline=center] (f){$\nu_{\bb A, \bb A}^{-1}$};
\draw ([xshift=-0.3cm]f.north) -- node[pos=0.5, left]{$\bbover A$}++(0, 0.75);
\draw ([xshift=0.3cm]f.north) -- node[pos=0.5, right]{$\bbover A$}++(0, 0.75);
\draw ([xshift=-0.55cm, yshift=1cm]f.north) node[draw=black,fill=white,very thick,rounded corners,baseline=center] (g){$\sigma_{\bb A}^{-1}$};
\draw ([xshift=0.55cm, yshift=1cm]f.north) node[draw=black,fill=white,very thick,rounded corners,baseline=center] (h){$\sigma_{\bb A}^{-1}$};
\filldraw ([xshift=0.55cm, yshift=0.5cm]g.north) circle (0pt) coordinate(i);
\draw (g.north) to [out=90, in=180] node[pos=0.7, above left]{$\bb A$} (i);
\draw (h.north) to [out=90, in=0] node[pos=0.7, above right]{$\bb A$} (i);
\draw (i) --node[pos=1, above]{$\bb A$}++(0, 0.5);
\draw ([xshift=-0.05cm]j.south) --++(0, -0.75);
\draw ([xshift=0.05cm]j.south) --node[pos=0.55, black, right]{$\overline{\bbover A \otimes \bb A}$}++(0, -0.75);
\draw ([yshift=-1.45cm]j) node[draw=black,fill=white,very thick,rounded corners,baseline=center] (k){$\nu_{\bbover A, \bb A}$};
\draw ([xshift=-0.3cm, yshift=-1cm]k.south) --node[pos=0, black, below] {$\bbover A$}++(0, 1);
\draw ([xshift=0.3cm, yshift=-0.25cm]k.south) --++(0, 0.25);
\node[] at ([xshift=0.75cm, yshift=0cm]k.south) {$\overline{\bbover A}$};
\draw ([xshift=0.3cm, yshift=-0.5cm]k.south) node[draw=black,fill=white,very thick,rounded corners,baseline=center] (l){$\varphi_{\bb A}$};
\draw ([yshift=-0.25cm]l.south) --node[pos=0, black, below] {$\bb A$}++(0, 0.25);

\node[] (m) at ([xshift=1.75cm]j){$=$};

\draw ([xshift=1.5cm]m) node[draw=black,fill=white,very thick,rounded corners,baseline=center] (n){$\nu_{\bbover A, \bb A}$};
\draw ([xshift=-0.3cm]n.south) --node[pos=1, black, below] {$\bbover A$}++(0, -2.85);
\draw ([xshift=0.3cm]n.south) --node[pos=0.5, right]{$\bbover A$} ++(0, -0.75);
\draw ([xshift=0.3cm, yshift=-1cm]n.south) node[draw=black,fill=white,very thick,rounded corners,baseline=center] (o){$\sigma^{-1}_{\bb A}$};
\draw (o.south) --node[pos=0.5, black, right] {$\overline{\bbover A}$}++(0, -0.75);
\draw ([yshift=-1cm]o.south) node[draw=black,fill=white,very thick,rounded corners,baseline=center] (p){$\varphi_{\bb A}$};
\draw (p.south) --node[pos=1, below]{$\bb A$}++(0, -0.25);

\draw ([xshift=-0.05cm]n.north) -- ++(0, 0.75);
\draw ([xshift=0.05cm]n.north) --node[pos=0.5, right]{$\overline{\bb A \otimes \bb A}$} ++(0, 0.75);
\draw ([yshift=1.5cm]n) node[draw=black,fill=white,very thick,rounded corners,baseline=center] (t){$\nu_{\bb A, \bb A}^{-1}$};
\draw ([xshift=-0.3cm]t.north) -- node[pos=0.5, left]{$\bbover A$}++(0, 0.75);
\draw ([xshift=0.3cm]t.north) -- node[pos=0.5, right]{$\bbover A$}++(0, 0.75);
\draw ([xshift=-0.55cm, yshift=1cm]t.north) node[draw=black,fill=white,very thick,rounded corners,baseline=center] (q){$\sigma_{\bb A}^{-1}$};
\draw ([xshift=0.55cm, yshift=1cm]t.north) node[draw=black,fill=white,very thick,rounded corners,baseline=center] (r){$\sigma_{\bb A}^{-1}$};
\filldraw ([xshift=0.55cm, yshift=0.5cm]q.north) circle (0pt) coordinate(s);
\draw (q.north) to [out=90, in=180] node[pos=0.7, above left]{$\bb A$} (s);
\draw (r.north) to [out=90, in=0] node[pos=0.7, above right]{$\bb A$} (s);
\draw (s) --node[pos=1, above]{$\bb A$}++(0, 0.5);

\node[] (t) at ([xshift=1.75cm]n){$=$};

\draw ([xshift=1.25cm]t.north) node[draw=black,fill=white,very thick,rounded corners,baseline=center] (u){$\sigma_{\bb A}^{-1}$};
\draw ([xshift=1.1cm]u) node[draw=black,fill=white,very thick,rounded corners,baseline=center] (v){$\sigma_{\bb A}^{-1}$};
\filldraw ([xshift=0.55cm, yshift=0.5cm]u.north) circle (0pt) coordinate (w);
\draw (u.north) to [out=90, in=180] node[pos=0.7, above left]{$\bb A$} (w);
\draw (v.north) to [out=90, in=0] node[pos=0.7, above right]{$\bb A$} (w);
\draw (w) --node[pos=1, above]{$\bb A$}++(0, 0.75);
\draw ([yshift=-1.25cm]v) node[draw=black,fill=white,very thick,rounded corners,baseline=center] (x){$\sigma_{\bb A}^{-1}$};
\draw(v.south) to node[pos=0.5, right]{$\bbover A$} (x.north);
\draw ([yshift=-1.25cm]x) node[draw=black,fill=white,very thick,rounded corners,baseline=center] (y){$\varphi_{\bb A}$};
\draw(x.south) to node[pos=0.5, right]{$\overline{\bbover A}$} (y.north);
\draw(y.south) --node[pos=1, below]{$\bb A$}++(0, -0.5);
\draw(u.south) to node[pos=1, below]{$\bbover A$} ([xshift=-1.1cm,yshift=-0.5cm]y.south);

\node[] (z) at ([xshift=1.25cm]v){$=$};

\filldraw ([xshift=1.5cm]z) circle (0pt) coordinate (aa);
\draw (aa) --node[pos=1,above]{$\bb A$}++(0, 0.75);
\draw ([xshift=-0.6cm,yshift=-1.25cm]aa) node[draw=black,fill=white,very thick,rounded corners,baseline=center] (ab){$\sigma_{\bb A}^{-1}$};
\draw (ab.north) to [out=90, in=180] node[pos=0.7, above left]{$\bb A$} (aa);
\draw (ab.south) --node[pos=1, below]{$\bbover A$}++(0, -0.75) coordinate (ac);
\draw ([xshift=1.2cm]ab.north) to [out=90, in=0] (aa);
\draw ([xshift=1.2cm]ab.north) to node[pos=1, below]{$\bb A$} ([xshift=1.2cm]ac);
\end{tikzpicture}
.
\end{center}
The diagram in \ref{H3} taken with the inner product on $\bb A$ is manifestly positive, as it lives in a corner of the $\rmC^*$-algebra $\cM_{\bb A}(\overline c \otimes c \otimes \bb A, \overline c \otimes c \otimes \bb A)$.
Similarly, $\bb A(a)$ is a corner of a $\rmC^*$ algebra (in $\cM_{\bb A}$), hence complete.
\end{ex}

Without too much more effort, we obtain the following example.

\begin{ex}  \label{ex:cyclicmodulecat}
Let $\bb X$ be a right $\bb A$ Hilbert $\rmC^*$ module object.
Then $\bb c \otimes \bb X$ is a right $\bb A$ Hilbert $\rmC^*$ module object as well.
The inner product is:
\[
\begin{tikzpicture}[baseline={(current bounding box.center)}, thick]
\draw (0,0) node[draw=black,fill=white,very thick,rounded corners,baseline=center, inner xsep=0.2cm] (a){$\langle \cdot | \cdot \rangle_{\bf A}^{\bf X}$};
\draw (a.north) --node[pos=1, above] {$\bb A$} ++(0, 0.25);
\draw [red, very thick] ([xshift=-0.5cm]a.south) -- node[pos=0.4, left, black] {$\bbover X$} ++ (0, -0.75);
\draw ([xshift=-0.5cm, yshift=-1.25cm]a) node[draw=black,fill=white,very thick,rounded corners,baseline=center] (b){$\nu^{-1}$};
\draw ([xshift=-0.05cm]b.south) --node[pos=1, below] {$\overline{\bb c \otimes \bb X}$} ++ (0, -0.5);
\draw [red, very thick] ([xshift=0.05cm]b.south) --++ (0, -0.5) coordinate (c);
\draw [red, very thick] ([xshift=0.5cm]a.south) to node[pos=1, below, black] {$\bb c \otimes \bb X$} ([xshift=0.95cm]c) coordinate (d);
\draw ([xshift=0.1cm]b.north) to [out=90, in=180] node[pos=0.3, above] {$\bbover c$} ([yshift=-0.3cm]a.south) coordinate (e);
\draw ([xshift=0.9cm]b.north) coordinate (f) to [out=90, in=0] (e);
\draw (f) to ([xshift=-0.1cm]d);
\end{tikzpicture}
\]
Right-linearity is clear.
By involutive category coherence equalities, we write 
\begin{center}
\begin{tikzpicture}[baseline={(current bounding box.center)}, thick]
\draw ([xshift=2cm]a) node[draw=black,fill=white,very thick,rounded corners,baseline=center] (b){$\overline{\langle\cdot|\cdot\rangle_{\mathbf{A}}^{\bb c \otimes \bb X}}$};
\draw [red, very thick] ([xshift=-0.05cm]b.south) --++(0, -0.75);
\draw [red, very thick] ([xshift=0.05cm]b.south) --node[pos=0.55, black, right]{$\overline{\bbover X \otimes \bb X}$}++(0, -0.75);
\draw (b.north) --node[pos=0.5, black, right] {$\bbover A$}++(0,0.5);
\draw ([yshift=1.25cm]b) node[draw=black,fill=white,very thick,rounded corners,baseline=center] (c){$\sigma^{-1}_{\bb A}$};
\draw (c.north)--node[pos=1, black, above] {$\bb A$}++(0, 0.5);
\draw ([yshift=-1.45cm]b) node[draw=black,fill=white,very thick,rounded corners,baseline=center] (d){$\nu_{\bbover X, \bb X}$};
\draw [red, very thick] ([xshift=-0.3cm, yshift=-1cm]d.south) --node[pos=0, black, below] {$\bbover X$}++(0, 1);
\draw [red,very thick] ([xshift=0.3cm, yshift=-0.25cm]d.south) --++(0, 0.25);
\node[] at ([xshift=0.75cm, yshift=0cm]d.south) {$\overline{\bbover X}$};
\draw ([xshift=0.3cm, yshift=-0.5cm]d.south) node[draw=black,fill=white,very thick,rounded corners,baseline=center] (e){$\varphi_{\bb X}$};
\draw [red, very thick] ([yshift=-0.25cm]e.south) --node[pos=0, black, below] {$\bb X$}++(0, 0.25);

\node[] (m) at ([xshift=2.25cm]b) {$=$};

\draw ([xshift=2.5cm, yshift=1.5cm]m) node[draw=black,fill=white,very thick,rounded corners,baseline=center] (n){$\overline{1_{\bbover X} \otimes \ev \otimes 1_{\bb X}}$};
\draw [red, very thick] ([xshift=-0.05cm]n.north) --++ (0, 0.75);
\draw [red, very thick] ([xshift=0.05cm]n.north) -- node[black, pos=0.4, right] {$\overline{\bbover X \otimes \bb X}$} ++ (0, 0.75);
\draw ([yshift=1.5cm]n) node[draw=black,fill=white,very thick,rounded corners,baseline=center] (o){$\overline{\langle \cdot | \cdot \rangle_{\bb A}^{\bb X}}$};
\draw (o.north) --node[pos=0.4, right] {$\bbover A$} ++ (0, 0.75);
\draw ([yshift=1.5cm]o) node[draw=black,fill=white,very thick,rounded corners,baseline=center] (p){$\sigma^{-1}_{\bb A}$};
\draw (p.north) --node[pos=1, above] {$\bb A$} ++ (0, 0.5);
\draw [red, very thick] ([xshift=-0.13cm]n.south) --++ (0, -0.75);
\draw ([xshift=-0.05cm]n.south) --++ (0, -0.75);
\draw ([xshift=0.05cm]n.south) --++ (0, -0.75);
\draw [red, very thick] ([xshift=0.13cm]n.south) -- node[pos=0.5, right, black] {$\overline{\overline{\bb c \otimes \bb X} \otimes \bb c \otimes \bb X}$} ++ (0, -0.75);
\draw ([yshift=-1.5cm]n) node[draw=black,fill=white,very thick,rounded corners,baseline=center] (q){$\overline{\nu^{-1} \otimes 1_{\bb c \otimes \bb X}}$};
\draw ([xshift=-0.15cm]q.south) --++ (0,-0.75);
\draw [red, very thick] ([xshift=-0.05cm]q.south) --++ (0,-0.75);
\draw ([xshift=0.05cm]q.south) --++ (0,-0.75);
\draw [red, very thick] ([xshift=0.15cm]q.south) -- node[pos=0.4, right, black] {$\overline{\overline{\bb c \otimes \bb X} \otimes \bb c \otimes \bb X}$} ++ (0,-0.75);
\draw ([yshift=-1.25cm]q) node[draw=black,fill=white,very thick,rounded corners,baseline=center, inner xsep=0.75cm] (z){$\nu$};
\draw ([xshift=0.55cm]z.south) --++ (0, -0.75);
\draw [red, very thick] ([xshift=0.65cm]z.south) --node[pos=0.4, right, black] {$\overline{\overline{\bb c \otimes \bb X}}$} ++ (0, -0.75);
\draw ([xshift=-0.55cm]z.south) --++ (0, -0.75);
\draw [red, very thick] ([xshift=-0.65cm]z.south) --node[pos=0.4, left, black] {$\overline{\bb c \otimes \bb X}$} ++ (0, -0.75);
\draw ([xshift=-0.6cm, yshift=-1.25cm]z) node[draw=black,fill=white,very thick,rounded corners,baseline=center] (r){$\nu$};
\draw ([xshift=0.6cm, yshift=-1.25cm]z) node[draw=black,fill=white,very thick,rounded corners,baseline=center] (s){$\nu$};
\draw ([xshift=-0.05cm]s.south) --++ (0, -0.75);
\draw [red, very thick] ([xshift=0.05cm]s.south) --node[pos=0.4, right, black] {$\overline{\bbover X \otimes \bbover c}$} ++ (0, -0.75);
\draw ([xshift=0.05cm]r.south) --++ (0, -0.75);
\draw [red, very thick] ([xshift=-0.05cm]r.south) --node[pos=0.4, left, black] {$\bbover X \otimes \bbover c$} ++ (0, -0.75);
\draw ([yshift=-1.25cm]r) node[draw=black,fill=white,very thick,rounded corners,baseline=center] (t){$\nu^{-1}$};
\draw ([yshift=-1.25cm]s) node[draw=black,fill=white,very thick,rounded corners,baseline=center, inner xsep=0.3cm] (u){$\nu$};
\draw ([xshift=-0.3cm]u.south) --node[pos=0.4, left] {$\overline{\bbover c}$} ++ (0, -0.75);
\draw [red, very thick] ([xshift=0.3cm]u.south) --node[pos=0.4, right, black] {$\overline{\bbover X}$} ++ (0, -0.75);
\draw ([xshift=0.3cm, yshift=-1.25cm]u) node[draw=black,fill=white,very thick,rounded corners,baseline=center] (v){$\varphi$};
\draw ([xshift=-0.3cm, yshift=-1.25cm]u) node[draw=black,fill=white,very thick,rounded corners,baseline=center] (w){$\varphi$};
\draw (w.south) -- node[pos=1, below, black] {$\bb c$} ++(0, -0.5) coordinate (x);
\draw [red, very thick] (v.south) -- node[pos=1, below, black] {$\bb X$} ++(0, -0.5);
\draw ([xshift=-0.1cm]t.south) to ([xshift=-1cm]x);
\draw [red, very thick] (t.south) to node[pos=1, below, black] {$\overline{\bb c \otimes \bb X}$} ([xshift=-0.9cm]x);
\end{tikzpicture}
\end{center}
Applying $(H3)$ for $\langle \cdot | \cdot \rangle_{\bf A}^{\bf X}$, it suffices to show the below equality of diagrams, which one can quickly check using naturality of $\nu$, and Equations \ref{eq:assoc} and \ref{eq:decomp}.
\begin{center}
\begin{tikzpicture}[baseline={(current bounding box.center)}, thick]
\draw (0, 0) node[draw=black,fill=white,very thick,rounded corners,baseline=center, inner xsep=0.5cm] (a){$\nu$};
\draw [red, very thick] ([xshift=-0.5cm]a.south) -- node[pos=1, below, black] {$\bbover X$} ++ (0, -0.75) coordinate (b);
\draw [red, very thick] ([xshift=0.5cm]a.south) -- node[pos=1, below, black] {$\overline{\bbover X}$} ++ (0, -0.75) coordinate (c);
\draw ([xshift=0.3cm]b) to [out=90, in=180] node[pos=0, below] {$\bbover c$} ([yshift=-0.3cm]a.south) coordinate (e);
\draw ([xshift=-0.3cm]c) to [out=90, in=0] node[pos=0, below] {$\bb c$} (e);
\draw [red, very thick] ([xshift=-0.1cm]a.north) --++ (0, 0.75) coordinate (f);
\draw [red, very thick] ([xshift=0.1cm]a.north) --++ (0, 0.75);
\node[above] at ([xshift=0.1cm]f) {$\overline{\bbover X \otimes \bb X}$};

\node[] (g) at ([xshift=1.5cm]a) {$=$};

\draw ([xshift=2.5cm]g) node[draw=black,fill=white,very thick,rounded corners,baseline=center, inner xsep=0.75cm] (h){$\nu$};
\draw [red, very thick] ([xshift=-0.15cm]h.north) --++ (0, 0.75);
\draw ([xshift=-0.05cm]h.north) --++ (0, 0.75);
\draw ([xshift=0.05cm]h.north) --++ (0, 0.75);
\draw [red, very thick] ([xshift=0.15cm]h.north) --node[pos=0.4, right, black] {$\overline{\bbover X \otimes \bbover c \otimes \bb c \otimes \bb X}$}++ (0, 0.75);
\draw ([yshift=1.25cm]h) node[draw=black,fill=white,very thick,rounded corners,baseline=center] (i){$\overline{1_{\bbover X} \otimes \ev \otimes 1_{\bb X}}$};
\draw [red, very thick] ([xshift=-0.1cm]i.north) --++ (0, 0.5) coordinate (k);
\draw [red, very thick] ([xshift=0.1cm]i.north) --++ (0, 0.5);
\node[above] at ([xshift=0.1cm]k) {$\overline{\bbover X \otimes \bb X}$};
\draw ([xshift=-0.65cm]h.south) -- node[pos=0.4, left] {$\overline{\bb c \otimes \bb X}$} ++ (0, -0.75);
\draw [red, very thick] ([xshift=-0.55cm]h.south) --++ (0, -0.75);
\draw ([xshift=0.65cm]h.south) -- node[pos=0.4, right] {$\overline{\bbover c \otimes \bbover X}$} ++ (0, -0.75);
\draw [red, very thick] ([xshift=0.55cm]h.south) --++ (0, -0.75);
\draw ([xshift=-0.6cm, yshift=-1cm]h) node[draw=black,fill=white,very thick,rounded corners,baseline=center, inner xsep=0.3cm] (l){$\nu$};
\draw [red, very thick] ([xshift=-0.2cm]l.south) --node[pos=1, below, black] {$\bbover X$} ++ (0, -1);
\draw ([xshift=0.2cm]l.south) --node[pos=1, below, black] {$\bbover c$} ++ (0, -1);
\draw ([xshift=0.6cm, yshift=-1cm]h) node[draw=black,fill=white,very thick,rounded corners,baseline=center, inner xsep=0.3cm] (m){$\nu$};
\draw ([xshift=-0.2cm]m.south) --++ (0, -0.5);
\node[left] at ([xshift=-0.3cm, yshift=-0.05cm]m.south){$\overline{\bbover c}$};
\draw ([xshift=-0.2cm, yshift=-0.75cm]m) node[draw=black,fill=white,very thick,rounded corners,baseline=center] (n){$\varphi$};
\draw (n.south) --node[pos=1, below] {$\bb c$} ++ (0, -0.25);
\draw [red, very thick]  ([xshift=0.2cm]m.south) --node[pos=1, below, black] {$\overline{\bbover X}$} ++ (0, -1);

\end{tikzpicture}
\end{center}
By a similar justification to Example \ref{ex:overself}, \ref{H3} and \ref{H4} are satisfied.
\end{ex}

Having defined Hilbert $\rmC^*$ module objects, we define their morphisms.

\begin{defn}
A right $\bb A$-linear morphism $\alpha$ is an adjointable operator if there exists an adjoint $\alpha^* \colon \bb Y_{\bb A} \to \bb X_{\bb A}$ satisfying the diagram
\begin{equation}
\label{eq:AdjOp}
\begin{tikzpicture}[baseline={(current bounding box.center)}, thick]
\draw (0,0) node[draw=black, fill=white, very thick, rounded corners, baseline=center] (a){$\langle \cdot | \cdot \rangle_{\bb A}^{\bb Y}$};
\draw ([xshift=-0.3cm, yshift=-1.25cm]a) node[draw=black, fill=white, very thick, rounded corners, baseline=center] (b){$\overline{\alpha^*}$};
\draw (a.north) --node[pos=1, above]{$\bb A$} ++(0, 0.5);
\draw [red, very thick] ([xshift=-0.3cm]a.south) to node[pos=0.4, left, black]{$\bbover Y$} (b.north);
\draw [red, very thick] (b.south) --node[pos=1, below, black]{$\bbover X$}++(0, -0.5) coordinate (c);
\draw [red, very thick] ([xshift=0.3cm]a.south) to node[pos=1, below, black]{$\bb Y$} ([xshift=0.6cm]c);

\node[] (d) at ([xshift=1cm]a) {$=$};

\draw ([xshift=1cm]d) node[draw=black, fill=white, very thick, rounded corners, baseline=center] (e){$\langle \cdot | \cdot \rangle_{\bb A}^{\bb X}$};
\draw ([xshift=0.3cm, yshift=-1.25cm]e) node[draw=black, fill=white, very thick, rounded corners, baseline=center] (f){$\alpha$};
\draw (e.north) --node[pos=1, above]{$\bb A$} ++(0, 0.5);
\draw [red, very thick] ([xshift=0.3cm]e.south) to node[pos=0.4, right, black]{$\bb X$} (f.north);
\draw [red, very thick] (f.south) --node[pos=1, below, black]{$\bb Y$}++(0, -0.5) coordinate (g);
\draw [red, very thick] ([xshift=-0.3cm]e.south) to node[pos=1, below, black]{$\bbover X$} ([xshift=-0.6cm]g);
\end{tikzpicture}
\end{equation}
such that both $\alpha$ and $\alpha^*$ are uniformly bounded, i.e. have finite norm $\|\alpha\| = \sup \limits_{c \in \cC} \|\alpha_c\|$, where the norm on $\alpha_c$ is the operator norm with respect to the norms on $\bb X(a)$ and $\bb Y(a)$.
\end{defn}

We make the following observations about this definition.

\begin{rmk}
For an $\bb A$ module morphism $\alpha$ satisfying the diagram \ref{eq:AdjOp}, $$\|a\| = \sup \limits_{c \in \cC} \|\alpha_c\| = \sup \limits_{c \in \Irr(\cC)} \|\alpha_c\|$$
This follows from semisimplicity of $\cC$.
If the $\bb A$ module morphism $\alpha$ satisfies diagram \ref{eq:AdjOp}, then $\|\alpha\| < \infty$ iff $\|\alpha^*\| < \infty$.
Similarly, the closed graph theory implies that an $\bb A$ module morphism satisfying \ref{eq:AdjOp} has bounded components.

Warning!
That an $\bb A$ module morphism $\alpha$ satisfies the diagram \ref{eq:AdjOp} does not imply that $\alpha$ is an adjointable operator.
Consider $\cC = \Hilb(\mathbb Z)$, $\bb X$ such that $\bb X (\mathbb{C}_n) = \mathbb{C}$, and $\bb A = \cC(1_\cC)$.
Then $\alpha \colon \bb X \to \bb X$ given by $\alpha_n \colon \bb X(\mathbb{C}_n) \to \bb X(\mathbb{C}_n)$ is multiplication by $n$ satisfies diagram \ref{eq:AdjOp}.
However $\alpha$ is not uniformly bounded as $\|\alpha_n\| = n$.
\end{rmk}

As module objects for algebra objects form a module category, so do Hilbert $\rmC^*$ module objects for a $\rmC^*$ algebra object.

\begin{prop} \label{prop:modulecat}
For a $\rmC^*$-algebra object $\bb A$, the category $\CMod_{\bb A}$ of right (resp. left) $\bb A$ Hilbert $\rmC^*$-modules with adjointable operators is a $\cC$ module category.
\end{prop}

\begin{proof}
We claim that $\otimes_{\CMod_{\bb A}} = \otimes_{\Vect(\cC)}$ is a choice of action of $\cC$ on $\CMod_{\bb A}$ which gives the result.
By Example \ref{ex:cyclicmodulecat}, $\bb c \otimes \bb X \in \CMod_{\bb A}$.
We show that for $f \colon c \to d$ in $\cC$ and an adjointable operator $\alpha \colon \bb X \to \bb Y$, $\bb f \otimes \alpha$ is an adjointable operator.
Its adjoint is $\bb f^* \otimes \alpha^*$, where $\bb f^*$ is precomposition by the adjoint (in $\cC$) of $f$.
By semisimplicity of $\cC$, $\bb f \otimes \alpha$ and $\bb f^* \otimes \alpha^*$ are uniformly bounded, and as $\alpha$ is right $\bb A$-linear, so too is $\bb f \otimes \alpha$.
To see that $\bb f^* \otimes \alpha^* = (\bb f \otimes \alpha)^*$, we compute

\begin{center}
\begin{tikzpicture}[baseline={(current bounding box.center)}, thick]
\draw (0,0) node[draw=black, fill=white, very thick, rounded corners, baseline=center, inner xsep=0.3cm] (a){$\langle \cdot | \cdot \rangle_{\bb A}$};
\draw [red, very thick] ([xshift=-0.6cm]a.south) --node[pos=0.3, left, black] {$\bbover X$} ++ (0, -0.75);
\draw ([xshift=-0.45cm, yshift=-0.65cm]a.south) to [out=90, in=180] node[pos=0.8, below] {$\bbover c$} ([yshift=-0.1cm]a.south) coordinate (c);
\draw ([xshift=0.45cm, yshift=-0.65cm]a.south) coordinate (f) to [out=90, in=0] (c);
\draw ([xshift=-0.6cm, yshift=-1.25cm]a) node[draw=black, fill=white, very thick, rounded corners, baseline=center] (b){$\nu^{-1}$};
\draw ([xshift=-0.05cm]b.south) --node[left, pos=0.4] {$\overline{\bb c \otimes \bb X}$} ++(0, -0.75);
\draw [red, very thick] ([xshift=0.05cm]b.south) --++(0, -0.75);
\draw ([yshift=-1.25cm]b) node[draw=black, fill=white, very thick, rounded corners, baseline=center] (d){$\overline{\bb f^* \otimes \alpha^*}$};
\draw ([xshift=-0.05cm]d.south) --++ (0, -0.5);
\draw [red, very thick] ([xshift=0.05cm]d.south) --++ (0, -0.5) coordinate (e);
\node[below] at ([yshift=-0.5cm]d.south) {$\overline{\bb d \otimes \bb Y}$};
\draw [red, very thick] ([xshift=0.6cm]a.south) to ([xshift=1.15cm]e);
\draw (f) to ([xshift=1cm]e);
\node[below] at ([xshift=1.1cm]e) {$\bb c \otimes \bb X$};
\draw (a.north) --node[pos=1, above] {$\bb A$}++ (0, 0.5);

\node[] (g) at ([xshift=1.5cm]a) {$=$};

\draw ([xshift=1.75cm]g) node[draw=black, fill=white, very thick, rounded corners, baseline=center, inner xsep=0.6cm] (a){$\langle \cdot | \cdot \rangle_{\bb A}$};
\draw [red, very thick] ([xshift=-0.9cm]a.south) --node[pos=0.3, left, black] {$\bbover X$} ++ (0, -0.75);
\draw ([xshift=-0.25cm, yshift=-1.25cm]a) node[draw=black, fill=white, very thick, rounded corners, baseline=center] (h){$\overline{\bb f^*}$};
\draw (h.south) --node[pos=0.4, right] {$\bbover d$} ++ (0, -0.75);
\draw ([xshift=-0.25cm, yshift=-0.65cm]a.south) to [out=90, in=180] node[pos=0.8, below] {$\bbover c$} ([xshift=0.2cm, yshift=-0.1cm]a.south) coordinate (c);
\draw ([xshift=0.65cm, yshift=-0.65cm]a.south) coordinate (f) to [out=90, in=0] (c);
\draw ([xshift=-0.95cm, yshift=-1.25cm]a) node[draw=black, fill=white, very thick, rounded corners, baseline=center] (b){$\overline{\alpha^*}$};
\draw [red, very thick] (b.south) --node[pos=0.4, left, black] {$\bbover Y$} ++(0, -0.75);

\draw ([xshift=0.325cm, yshift=-1.25cm]b) node[draw=black, fill=white, very thick, rounded corners, baseline=center, inner xsep=0.3cm] (d){$\nu^{-1}$};
\draw ([xshift=-0.05cm]d.south) --++ (0, -0.5);
\draw [red, very thick] ([xshift=0.05cm]d.south) --++ (0, -0.5) coordinate (e);
\node[below] at ([yshift=-0.5cm]d.south) {$\overline{\bb d \otimes \bb Y}$};
\draw [red, very thick] ([xshift=0.9cm]a.south) to ([xshift=1.5cm]e);
\draw (f) to ([xshift=1.25cm]e);
\node[below] at ([xshift=1.4cm]e) {$\bb c \otimes \bb X$};
\draw (a.north) --node[pos=1, above] {$\bb A$}++ (0, 0.5);
\end{tikzpicture}
\end{center}

That $\bb f \otimes \alpha$ is adjointable follows from an application of diagram \ref{eq:AdjOp} for $\alpha$ and the compatibility of the dagger and rigid structures of $\cC$.
Bifunctoriality of $\otimes_{\CMod_{\bb A}}$ follows from here.
The associator is given by the associator for $\Vect(\cC)$, hence satisfies the pentagon axiom.
Its adjoint is its inverse, as seen in the following computation (suppressing the involutive structure).
\begin{center}
\begin{tikzpicture}[baseline={(current bounding box.center)}, thick]
\draw (0, 0) node[draw=black,fill=white,very thick,rounded corners,baseline=center] (a){$\overline \xi$};
\filldraw [black] ([xshift=0.45cm]a.south) circle (0pt) coordinate (b);
\filldraw [black] ([xshift=0.85cm]a.south) circle (0pt) coordinate (c);
\filldraw [black] ([xshift=0.65cm, yshift=-0.5cm]a.south) circle (0pt) coordinate (d);
\node[above] at (d) {$\gamma$};
\filldraw [black] ([xshift=0.35cm, yshift=-1cm]a.south) circle (0pt) coordinate (e);
\node[above] at (e) {$\beta$};
\draw (a.south) to [out=-90, in=180] node[pos=0.4, left] {$\bbover b$} (e);
\draw (b) to [out=-90, in=180] node[pos=0.4, left] {$\bbover c$} (d);
\draw (c) to [out=-90, in=0] node[pos=0.4, right] {$\bbover d$} (d);
\draw (d) to [out=-90, in=0] node[pos=0.8, right] {$\bbover f$} (e);
\draw (e) --node[pos=1, below] {$\bbover a$}++ (0, -0.5) coordinate (j);
\draw (a.north) --node[pos=1, above] {$\bbover X$}++ (0, 0.5) coordinate (k);

\draw ([xshift=2cm, yshift=-0.25cm]a) node[draw=black,fill=white,very thick,rounded corners,baseline=center, inner xsep =0.5cm] (f){$\alpha$};
\draw ([xshift=0.35cm, yshift=-0.75cm]f) node[draw=black,fill=white,very thick,rounded corners,baseline=center] (g){$\eta$};
\filldraw [black] ([xshift=0cm, yshift=-1.25cm]f.south) circle (0pt) coordinate (h);
\node[below left] at (h) {$\delta$};
\filldraw [black] ([xshift=0.2cm, yshift=-1cm]f.south) circle (0pt) coordinate (i);
\node[right] at (i) {$\varepsilon$};
\draw (h) to [out=0, in=-90] node[pos=0.7, below right] {$\bb h$} (i);
\draw (h) to [out=180, in=-90] node[pos=0.6, left] {$\bb d$} ([xshift=-0.35cm]f.south);
\draw (i) to [out=0, in=-90] (g.south);
\draw (i) to [out=180, in=-90] ([xshift=-0.05cm]f.south);
\draw (g.north) --++ (0, 0.25);
\draw (h) to node[pos=0.9, below] {$\bb e$} ([xshift=1.65cm, yshift=-0.25cm]j);

\draw ([xshift=0.35cm]f.north) to node[pos=1, above] {$\bb X$} ([xshift=2.38cm]k);
\filldraw [black] ([yshift=0.7cm, xshift=0.45cm]c) circle (0pt) coordinate (l);
\filldraw [black] ([xshift=0.45cm, yshift=1cm]c) circle (0pt) coordinate (m);
\draw (c) to [out=90, in=180] (l);
\draw (b) to [out=90, in=180] (m);
\draw ([xshift=-0.05cm]f.north) to [out=90, in=0] node[pos=0.7, above] {$\bb c$} (m);
\draw ([xshift=-0.35cm]f.north) to [out=90, in=0] (l);

\node[] (m) at ([xshift=1.25cm]f.east) {$= \sum \limits_{i \in \Irr(\cC)} U_{\delta, \varepsilon}^{\zeta, \theta}$};

\draw ([xshift=1.5cm, yshift=0.25cm]m) node[draw=black,fill=white,very thick,rounded corners,baseline=center] (a){$\overline \xi$};
\filldraw [black] ([xshift=0.45cm]a.south) circle (0pt) coordinate (b);
\filldraw [black] ([xshift=0.85cm]a.south) circle (0pt) coordinate (c);
\filldraw [black] ([xshift=0.65cm, yshift=-0.5cm]a.south) circle (0pt) coordinate (d);
\node[above] at (d) {$\gamma$};
\filldraw [black] ([xshift=0.35cm, yshift=-1cm]a.south) circle (0pt) coordinate (e);
\node[above] at (e) {$\beta$};
\draw (a.south) to [out=-90, in=180] node[pos=0.4, left] {$\bbover b$} (e);
\draw (b) to [out=-90, in=180] node[pos=0.4, left] {$\bbover c$} (d);
\draw (c) to [out=-90, in=0] node[pos=0.4, right] {$\bbover d$} (d);
\draw (d) to [out=-90, in=0] node[pos=0.8, right] {$\bbover f$} (e);
\draw (e) --node[pos=1, below] {$\bbover a$}++ (0, -0.5) coordinate (y);
\draw (a.north) --node[pos=1, above] {$\bbover X$}++ (0, 0.5) coordinate (z);

\draw ([xshift=2.75cm, yshift=-0.1cm]a) node[draw=black,fill=white,very thick,rounded corners,baseline=center] (f){$\eta$};
\filldraw [black] ([xshift=-0.45cm]f.south) circle (0pt) coordinate (g);
\filldraw [black] ([xshift=-0.85cm]f.south) circle (0pt) coordinate (h);
\filldraw [black] ([xshift=-0.65cm, yshift=-0.5cm]f.south) circle (0pt) coordinate (i);
\node[above] at (i) {$\zeta$};
\filldraw [black] ([xshift=-0.35cm, yshift=-1cm]f.south) circle (0pt) coordinate (j);
\node[above] at (j) {$\theta$};
\draw (f.south) to [out=-90, in=0] node[pos=0.4, right] {$\bb g$} (j);
\draw (g) to [out=-90, in=0] node[pos=0.4, right] {$\bb c$} (i);
\draw (h) to [out=-90, in=180] node[pos=0.4, left] {$\bb d$} (i);
\draw (i) to [out=-90, in=180] node[pos=0.8, left] {$\bb i$} (j);
\draw (j) to node[pos=1, below] {$\bb e$} ([xshift=2.05cm]y);
\draw (f.north) to node[pos=1, above] {$\bb X$} ([xshift=2.75cm]z);

\filldraw [black] ([xshift=0.55cm, yshift=0.5cm]c) circle (0pt) coordinate (k);
\filldraw [black] ([xshift=0.55cm, yshift=0.85cm]c) circle (0pt) coordinate (l);
\draw (b) to [out=90, in=180] (l);
\draw (c) to [out=90, in=180] (k);
\draw (g) to [out=90, in=0] (l);
\draw (h) to [out=90, in=0] (k);

\node[] (m) at ([xshift=0.75cm]f) {$=$};

\draw ([xshift=0.75cm, yshift=0.25cm]m) node[draw=black,fill=white,very thick,rounded corners,baseline=center] (a){$\overline \xi$};
\filldraw [black] ([xshift=0.45cm]a.south) circle (0pt) coordinate (b);
\filldraw [black] ([xshift=0.85cm]a.south) circle (0pt) coordinate (c);
\filldraw [black] ([xshift=0.65cm, yshift=-0.5cm]a.south) circle (0pt) coordinate (d);
\node[above] at (d) {$\gamma$};
\filldraw [black] ([xshift=0.35cm, yshift=-1cm]a.south) circle (0pt) coordinate (e);
\node[above] at (e) {$\beta$};
\draw (a.south) to [out=-90, in=180] node[pos=0.4, left] {$\bbover b$} (e);
\draw (b) to [out=-90, in=180] node[pos=0.4, left] {$\bbover c$} (d);
\draw (c) to [out=-90, in=0] node[pos=0.4, right] {$\bbover d$} (d);
\draw (d) to [out=-90, in=0] node[pos=0.8, right] {$\bbover f$} (e);
\draw (e) --node[pos=1, below] {$\bbover a$}++ (0, -0.5) coordinate (y);
\draw (a.north) --node[pos=1, above] {$\bbover X$}++ (0, 0.5) coordinate (z);

\draw ([xshift=2.75cm, yshift=-0.1cm]a) node[draw=black,fill=white,very thick,rounded corners,baseline=center] (f){$\eta$};
\filldraw [black] ([xshift=-0.45cm]f.south) circle (0pt) coordinate (g);
\filldraw [black] ([xshift=-0.85cm]f.south) circle (0pt) coordinate (h);
\filldraw [black] ([xshift=-0.25cm, yshift=-0.5cm]f.south) circle (0pt) coordinate (i);
\node[above] at (i) {$\varepsilon$};
\filldraw [black] ([xshift=-0.5cm, yshift=-1cm]f.south) circle (0pt) coordinate (j);
\node[above] at (j) {$\delta$};
\draw (g) to [out=-90, in=180] node[pos=0.4, left] {$\bb c$} (i);
\draw (f.south) to [out=-90, in=0] node[pos=0.4, right] {$\bb g$} (i);
\draw (j) to [out=0, in=-90] node[pos=0.6, right] {$\bb h$} (i);
\draw (j) to [out=180, in=-90] node[pos=0.7, left] {$\bb d$} (h);
\draw (j) to node[pos=1, below] {$\bb e$} ([xshift=1.9cm]y);
\draw (f.north) to node[pos=1, above] {$\bb X$} ([xshift=2.75cm]z);

\filldraw [black] ([xshift=0.55cm, yshift=0.5cm]c) circle (0pt) coordinate (k);
\filldraw [black] ([xshift=0.55cm, yshift=0.85cm]c) circle (0pt) coordinate (l);
\draw (b) to [out=90, in=180] (l);
\draw (c) to [out=90, in=180] (k);
\draw (g) to [out=90, in=0] (l);
\draw (h) to [out=90, in=0] (k);
\end{tikzpicture}

\begin{tikzpicture}[baseline={(current bounding box.center)}, thick]
\node[] (m) at (0,0) {$= \sum \limits_{i \in \Irr(\cC)} U_{\beta, \gamma}^{\kappa, \lambda}$};

\draw ([xshift=1.75cm, yshift=0.25cm]m) node[draw=black,fill=white,very thick,rounded corners,baseline=center] (a){$\overline \xi$};
\filldraw [black] ([xshift=0.45cm]a.south) circle (0pt) coordinate (b);
\filldraw [black] ([xshift=0.9cm]a.south) circle (0pt) coordinate (c);
\filldraw [black] ([xshift=0.25cm, yshift=-0.5cm]a.south) circle (0pt) coordinate (d);
\node[above] at (d) {$\kappa$};
\filldraw [black] ([xshift=0.5cm, yshift=-1cm]a.south) circle (0pt) coordinate (e);
\node[above] at (e) {$\lambda$};
\draw (a.south) to [out=-90, in=180] node[pos=0.4, left] {$\bbover b$} (d);
\draw (b) to [out=-90, in=0] node[pos=0.4, right] {$\bbover c$} (d);
\draw (c) to [out=-100, in=0] node[pos=0.4, right] {$\bbover d$} (e);
\draw (d) to [out=-90, in=180] node[pos=0.8, left] {$\bbover i$} (e);
\draw (e) --node[pos=1, below] {$\bbover a$}++ (0, -0.5) coordinate (y);
\draw (a.north) --node[pos=1, above] {$\bbover X$}++ (0, 0.5) coordinate (z);

\draw ([xshift=2.75cm, yshift=-0.1cm]a) node[draw=black,fill=white,very thick,rounded corners,baseline=center] (f){$\eta$};
\filldraw [black] ([xshift=-0.45cm]f.south) circle (0pt) coordinate (g);
\filldraw [black] ([xshift=-0.9cm]f.south) circle (0pt) coordinate (h);
\filldraw [black] ([xshift=-0.25cm, yshift=-0.5cm]f.south) circle (0pt) coordinate (i);
\node[above] at (i) {$\varepsilon$};
\filldraw [black] ([xshift=-0.5cm, yshift=-1cm]f.south) circle (0pt) coordinate (j);
\node[above] at (j) {$\delta$};
\draw (g) to [out=-90, in=180] node[pos=0.4, left] {$\bb c$} (i);
\draw (f.south) to [out=-90, in=0] node[pos=0.4, right] {$\bb g$} (i);
\draw (j) to [out=0, in=-90] node[pos=0.6, right] {$\bb h$} (i);
\draw (j) to [out=180, in=-80] node[pos=0.7, left] {$\bb d$} (h);
\draw (j) to node[pos=1, below] {$\bb e$} ([xshift=1.75cm]y);
\draw (f.north) to node[pos=1, above] {$\bb X$} ([xshift=2.75cm]z);

\filldraw [black] ([xshift=0.5cm, yshift=0.5cm]c) circle (0pt) coordinate (k);
\filldraw [black] ([xshift=0.5cm, yshift=0.85cm]c) circle (0pt) coordinate (l);
\draw (b) to [out=90, in=180] (l);
\draw (c) to [out=80, in=180] (k);
\draw (g) to [out=90, in=0] (l);
\draw (h) to [out=100, in=0] (k);

\node[] (m) at ([xshift=0.5cm, yshift=-0.25cm ]f.east) {$=$};

\draw ([xshift=2cm]f) node[draw=black,fill=white,very thick,rounded corners,baseline=center, inner xsep =0.5cm] (a){$\alpha^{-1}$};
\draw ([xshift=-0.35cm, yshift=-0.9cm]a) node[draw=black,fill=white,very thick,rounded corners,baseline=center] (b){$\overline \xi$};
\filldraw [black] ([xshift=0.25cm, yshift=-0.85cm]a.south) circle (0pt) coordinate (c);
\node[above] at (c) {$\gamma$};
\filldraw [black] ([xshift=-0.1cm, yshift=-1.35cm]a.south) circle (0pt) coordinate (d);
\node[above] at (d) {$\beta$};
\draw (d) to [out=0, in=-90] node[pos=0.6, right] {$\bbover f$} (c);
\draw (c) to [out=0, in=-90] node[pos=0.6, right] {$\bbover d$} ([xshift=0.5cm]a.south);
\draw (d) to [out=180, in=-90] node[pos=0.4, left] {$\bbover b$} (b.south);
\draw (c) to [out=180, in=-90] node[pos=0.85, right] {$\bbover c$} ([xshift=0cm]a.south);
\draw (b.north) --++ (0, 0.25);
\draw (d) --node[pos=0.9, below] {$\bbover a$}++ (0, -0.5) coordinate (m);

\draw ([xshift=2.75cm, yshift=0cm]a) node[draw=black,fill=white,very thick,rounded corners,baseline=center] (f){$\eta$};
\filldraw [black] ([xshift=-0.7cm]f.north) circle (0pt) coordinate (g);
\filldraw [black] ([xshift=-1.2cm]f.north) circle (0pt) coordinate (h);
\filldraw [black] ([xshift=-0.3cm, yshift=-0.5cm]f.south) circle (0pt) coordinate (i);
\node[above] at (i) {$\varepsilon$};
\filldraw [black] ([xshift=-0.65cm, yshift=-1cm]f.south) circle (0pt) coordinate (j);
\node[above] at (j) {$\delta$};
\draw (j) to [out=0, in=-90] node[pos=0.4, right] {$\bb h$} (i);
\draw (j) to [out=180, in=-90] node[pos=0.4, left] {$\bb d$} (h);
\draw (i) to [out=0, in=-90] node[pos=0.4, right] {$\bb g$} (f.south);
\draw (i) to [out=180, in=-90] node[pos=0.4, left] {$\bb d$} (g);
\draw (j) to node[pos=1, below] {$\bb e$} ([xshift=2.2cm]m);

\filldraw [black] ([xshift=1cm, yshift=0.35cm]a.north) circle (0pt) coordinate (k);
\filldraw [black] ([xshift=1cm, yshift=0.75cm]a.north) circle (0pt) coordinate (l);
\draw ([xshift=0.5cm]a.north) to [out=90, in=180] (k);
\draw ([xshift=0cm]a.north) to [out=90, in=180] (l);
\draw (h) to [out=90, in=0] (k);
\draw (g) to [out=90, in=0] (l);

\draw ([xshift=-0.35cm]a.north) --node[pos=1, above] {$\bbover X$}++ (0, 0.75) coordinate (n);
\draw (f.north) to node[pos=1, above] {$\bb X$} ([xshift=3.1cm]n);

\end{tikzpicture}
\end{center}
The unitor is given by the unitor for $\Vect(\cC)$, hence satisfies the triangle axiom.
Its adjoint is its inverse.
\end{proof}

In the theory of honest Hilbert $\rmC^*$ modules, the adjointable operators of a module form a $\rmC^*$ algebra.
Correspondingly, we would like the adjointable operators for a Hilbert $\rmC^*$ module object to form a $\rmC^*$ algebra object.
To prove this, we begin with a lemma.

\begin{lemma} \label{lem:directsums} 
The category $\CMod_{\bb A}$ admits direct sums.
\end{lemma}

\begin{proof}
For $\bb X, \bb Y \in \Vect(\cC)$, $\bb X \oplus \bb Y(c) \coloneq \bb X(c) \oplus \bb Y(c)$.
For $\bb X, \bb Y \in \CMod_{\bb A}$, we define a right $\bb A$-valued inner product on $\bb X \oplus \bb Y$ by $\langle \cdot | \cdot \rangle_{\bb A}^{\bb X \oplus \bb Y} = \langle \cdot | \cdot \rangle_{\bb A}^{\bb X} + \langle \cdot | \cdot \rangle_{\bb A}^{\bb Y}$.
One can quickly see that this is in fact a right $\bb A$-valued inner product.
Now suppose that there exist adjointable operators $\beta_{\bb X}$ and $\beta_{\bb Y}$ from $\bb X$ and $\bb Y$, respectively, to $\bb Z$.
The inclusions $\alpha_{\bb X}$ and $\alpha_{\bb Y}$ of $\bb X$ and $\bb Y$, respectively, into $\bb X \oplus \bb Y$ are adjointable operators.
We claim that $\gamma \coloneq \beta_{\bb X} \oplus \beta_{\bb Y}$ satisfies the universal property of the direct sum.
Since $\beta_{\bb X}$ and $\beta_{\bb Y}$ are adjointable operators, so is $\gamma$.
For $\bb X$, and similarly for $\bb Y$, $\gamma \circ \alpha_{\bb X} = \beta_{\bb X}$.
Since the construction in $\CMod_{\bb A}$ lifts the construction in $\Vect(\cC)$, $\gamma$ is the unique map satisfying the direct sum diagram.
\end{proof}

\begin{prop} \label{prop:adjopcstar}
The category $\CMod_{\bb A}$ is a $\rmC^*$ category.
\end{prop}

\begin{proof}
By Lemma \ref{lem:directsums}, $\CMod_{\bb A}$ admits direct sums.
So it suffices to show $B = \CMod_{\bb A} (\bb X, \bb X)$ is a $\rmC^*$-algebra for all $\bb X$ in $\CMod_{\bb A}$.
Clearly $B$ is a *-algebra.
Suppose $\sum_n \| \alpha_n\| < \infty$ for a sequence $(\alpha_n) \subset B$.
Then $\sum \|\alpha_{n,c}\| < \infty$ for all $c \in \cC$.
As $\bb X(c)$ is complete with respect to $\|\cdot \|_{\bb A(1_{\cC})}$, so is $\Vect(\bb X(c), \bb X(c))$.
So $\sum \alpha_{n,c}$ converges, say to $\alpha_c$.
Thus we define $\alpha \coloneq \sum \alpha_n$.
By linearity, $\alpha$ is a natural transformation.
That is, for $f \colon c \to d$, $$f\alpha_c = f\sum \alpha_{n, c} = \sum f \alpha_{n,c} = \sum \alpha_{n,d}f = (\sum \alpha_{n,d})f = \alpha_d f$$
One can easily see that $\alpha$ is right $\bb A$-linear.
We can find $\alpha^*$ in a similar way, and this is in fact the adjoint of $\alpha$ by continuity of the inner product.
More concretely, for $x, y \in \bb X(c)$.
$$\langle \alpha_c x | y \rangle = \left \langle \sum \alpha_{n,c} x | y \right \rangle = \sum \langle \alpha_{n,c} x | y \rangle = \sum \langle x | \alpha_{n,c}^* y \rangle = \langle x | \sum \alpha_{n,c}^* y \rangle = \langle x | \alpha^*_c y \rangle$$
Note that $\alpha$, and likewise $\alpha^*$, are both uniformly bounded.
Indeed, as $\sum \|\alpha_n\| < \infty$, $\|\alpha\| = \|\sum \alpha_n \| \leq \sum \|\alpha_n\| < \infty$.
Thus, $B$ is complete.
For $\alpha, \beta \in B$, one can quickly see that $\|\alpha \beta \| \leq \| \alpha \| \| \beta\|$.
Now we check $$\|\alpha_cx\|^2 = \langle \alpha_c x | \alpha_c x \rangle_{\bb A(1_{\cC})} = \langle \alpha^*_c \alpha_c x | x \rangle_{\bb A(1_{\cC})} \leq \|\alpha_C^*\alpha_c\| \|x\|^2$$
Taking supremumssupremums over $x \in \bb X(c)$ with $\|x\| \leq 1$, we see that $\|\alpha_c\|^2 \leq \|\alpha_c^*\alpha_c\|$ for all $c$.
Thus $\|\alpha\|^2 \leq \|\alpha^* \alpha\|$, and $B$ is a $\rmC^*$-algebra.
\end{proof}

With this in hand, we form the $\rmC^*$ algebra object of adjointable operators on a given Hilbert $\rmC^*$ module object.

\begin{defn} \label{defn:intend}
Let $\bb X_{\bb A}$ be a right Hilbert $\rmC^*$ module object.
By Proposition \ref{prop:modulecat}, the cyclic $\cC$ module category with basepoint $\bb X$, as a subcategory of $\CMod_{\bb A}$, is a $\cC$ module category.
By Proposition \ref{prop:adjopcstar}, this category is a $\rmC^*$ category.
Thus the corresponding algebra object is a $\rmC^*$ algebra object.
We call this category $\cM_{\bb X_{\bb A}}^{\rmC^*}$ and this object $\underline \End^{C^*}_{\cC}(\bb X_{\bb A})$.
Concretely, $\underline \End^{C^*}_{\cC}(\bb X_{\bb A})(c)$ is comprised of adjointable operators $\alpha \colon \bb c \otimes \bb X \to \bb X$.
\end{defn}

\begin{ex} \label{ex:algebraoverself}
By construction of $\underline \End^{C^*}_{\cC}$, we can see immediately that $\underline \End^{C^*}_{\cC}(\bb A_{\bb A}) \cong \bb A$.
In more detail, we can identify an adjointable operator $\alpha \colon \bb a \otimes \bb A \to \bb A$ with an element of $\bb A(a)$ under the isomorphism $\cM_{\bb A}(a \otimes \bb A, \bb A) \cong \bb A(a)$.
The elements $\alpha^* \colon \bb A \to \bb a \otimes \bb A$ and $j_a(\alpha)$ agree.
Effectively, this provides an isomorphism of $\cM_{\bfA}$ and $\cM_{\bfA}^{\rmC^*}$.
\end{ex}

Now we can define 1-morphisms in $\CAlg_{\cC}$.

\begin{defn} \label{def:HilbCBimodObj}
The right $\bb B$ Hilbert $\rmC^*$ module object $\bb X$ is an $\bb A-\bb B$ Hilbert $\rmC^*$ bimodule object when there exists a unital $*$-algebra natural transformation of $\rmC^*$ algebra objects $\phi^{\bb X} \colon \bb A \to \underline \End^{C^*}_{\cC}(\bb X_{\bb B})$.\footnote{By a unital algebra natural transformation $\alpha \colon \bf A \to \bf B$, we mean $\alpha(\iota_{\bfA}) = \iota_{\bfB}$. We make this assumption to avoid degeneracy issues.}
As usual, if the module object $\bb X$ is clear, we write $\phi^{\bb X}$ as $\phi$.
\end{defn}

\begin{ex}
A $\rmC^*$ module object $\bb A$ is an $\bb A-\bb A$ Hilbert $\rmC^*$ bimodule object.
\end{ex}

\begin{defn} 
Let $\bb X$ and $\bb Y$ be $\bb A-\bb B$ Hilbert $\rmC^*$ bimodule objects.
An adjointable operator $\alpha \colon \bb X \to \bb Y$ is a morphism of Hilbert $\rmC^*$ bimodule objects if
\begin{center}
\begin{tikzpicture}[baseline={(current bounding box.center)}, thick]
\draw (0,0) node[draw=black, fill=white, very thick, rounded corners, baseline=center] (a){$\alpha$};
\draw [red, very thick] (a.north) --node[pos=1, above, black] {$\bb Y$}++ (0, 0.75);
\draw [red, very thick] (a.south) -- node[pos=1, below, black] {$\bb X$} ++ (0, -2) coordinate (e);
\draw ([xshift=-0.75cm, yshift=-1cm]a.south) coordinate (c) to [out = 90, in = 180] node[pos=0.55, left] {$\underline \End^{C^*}_{\cC}(\bb X_{\bb B})$} ([yshift=-0.25cm]a.south);
\draw ([yshift=-0.25cm]c) node[draw=black, fill=white, very thick, rounded corners, baseline=center] (d){$\phi^{\bb X}$};
\draw (d.south) to node[pos=1, below] {$\bb A$} ([xshift=-0.75cm]e);

\node[] (f) at ([xshift=0.75cm, yshift=-0.25cm]a.south) {$=$};

\draw ([xshift=3cm]f) node[draw=black, fill=white, very thick, rounded corners, baseline=center] (a){$\alpha$};
\draw [red, very thick] (a.north) --node[pos=1, above, black] {$\bb Y$}++ (0, 1.25);
\draw [red, very thick] (a.south) -- node[pos=1, below, black] {$\bb X$} ++ (0, -1.5) coordinate (e);
\draw ([xshift=-0.75cm, yshift=-0.5cm]a.south) coordinate (c) to [out = 90, in = 180] node[pos=0.75, left] {$\underline \End^{C^*}_{\cC}(\bb Y_{\bb B})$} ([yshift=0.25cm]a.north);
\draw ([yshift=-0.25cm]c) node[draw=black, fill=white, very thick, rounded corners, baseline=center] (d){$\phi^{\bb Y}$};
\draw (d.south) to node[pos=1, below] {$\bb A$} ([xshift=-0.75cm]e);

\end{tikzpicture}
\end{center}
\end{defn}

To form a (weak) 2-category, we require composition of both 1-morphisms and 2-morphisms.
Composition of 2-morphisms is composition of natural transformations.
To compose 1-morphisms, we construct the relative tensor product of Hilbert $\rmC^*$ bimodule objects.
A technical tool for this construction is the completion of a pre-Hilbert modules object to a Hilbert module object.
See Appendix \ref{app:completionconstruction} for more details.

\begin{const} \label{const:relprod}
Let $\bb X$ be an $\bb A-\bb B$ Hilbert $\rmC^*$ bimodule object, let $\bb Y$ be a $\bb B-\bb C$ Hilbert $\rmC^*$ bimodule object, and let $\psi \colon \bb B \to \underline \End^{C^*}_{\cC}(\bb Y_{\bb C})$ be a *-algebra natural transformation.
The following form is a right $\bb C$-valued pre-inner product on $\bb X \otimes \bb Y$.

\begin{center}
\begin{tikzpicture}[baseline={(current bounding box.center)}, thick]
\draw (0,0) node[draw=black, fill=white, very thick, rounded corners, baseline=center, inner xsep=0.55cm] (a){$\langle \cdot | \cdot \rangle_{\bb C}$};
\draw [red, very thick] ([xshift=-0.75cm]a.south) -- node[pos=0.5, left, black] {$\bbover Y$}++ (0, -3.75);
\draw (a.north) -- node[pos=1, above, black] {$\bb C$} ++ (0, 0.5);
\draw [red, very thick] ([xshift=0.75cm]a.south) -- coordinate (b) node[pos=0.8, right, black] {$\underline \End^{C^*}_{\cC}(\bb Y_{\bb C}) \otimes \bb Y$}  ++ (0, -1);
\draw ([yshift=-1cm]a.south) coordinate (c) to [out=90, in=180] (b);
\draw ([yshift=-0.25cm]c) node[draw=black, fill=white, very thick, rounded corners, baseline=center] (d){$\psi$};
\draw (d.south) -- node[pos=0.4, right] {$\bb B$} ++ (0, -0.75);
\draw ([yshift=-1.25cm]d) node[draw=black, fill=white, very thick, rounded corners, baseline=center] (e){$\langle \cdot | \cdot \rangle_{\bb B}$};
\draw [red, very thick] ([xshift=-0.3cm]e.south) -- node[pos=0.4, right, black] {$\bbover X$} ++ (0, -0.75);
\draw ([xshift=-0.5cm, yshift=-1cm]e.south) node[draw=black, fill=white, very thick, rounded corners, baseline=center] (f){$\nu^{-1}$};
\draw [red, very thick] ([xshift=-0.05cm]f.south) --++ (0, -0.5);
\draw [red, very thick] ([xshift=0.05cm]f.south) --++ (0, -0.5) coordinate (g);
\node[below] at ([yshift=-0.5cm]f.south) {$\overline{\bb X \otimes \bb Y}$};
\draw [red, very thick] ([xshift=0.75cm, yshift=-1cm]a.south) to node[pos=1, below, black] {$\bb Y$} ([xshift=1.15cm]g);
\draw [red, very thick]  ([xshift=0.3cm]e.south) to node[pos=1, below, black] {$\bb X$} ([xshift=0.75cm]g);
\end{tikzpicture}
\end{center}

Right-linearity follows by right-linearity of $\langle \cdot | \cdot \rangle_{\bb C}$, and as $\underline \End^{C^*}_{\cC}(\bb Y_{\bb C})$ is comprised of right $\bb C$-module maps.
A brief computation shows anti-linearity.
It suffices to check positivity of the following diagram
\begin{center}
\begin{tikzpicture}[baseline={(current bounding box.center)}, thick]
\draw (0,0) node[draw=black, fill=white, very thick, rounded corners, baseline=center, inner xsep=0.55cm] (a){$\langle \cdot | \cdot \rangle_{\bb C}$};
\draw [red, very thick] ([xshift=-0.8cm]a.south) -- node[pos=0.5, left, black] {$\bbover Y$}++ (0, -3.75);
\draw (a.north) -- node[pos=1, above, black] {$\bb C$} ++ (0, 0.5);
\draw [red, very thick] ([xshift=0.8cm]a.south) -- coordinate (b) node[pos=0.8, right, black] {$\underline \End^{C^*}_{\cC}(\bb Y_{\bb C}) \otimes \bb Y$}  ++ (0, -1);
\draw ([yshift=-1cm]a.south) coordinate (c) to [out=90, in=180] (b);
\draw ([yshift=-0.25cm]c) node[draw=black, fill=white, very thick, rounded corners, baseline=center] (d){$\psi$};
\draw (d.south) -- node[pos=0.4, right] {$\bb B$} ++ (0, -0.75);
\draw ([yshift=-1.25cm]d) node[draw=black, fill=white, very thick, rounded corners, baseline=center] (e){$\langle \cdot | \cdot \rangle_{\bb B}$};
\draw [red, very thick] ([xshift=-0.25cm]e.south) -- node[pos=0.4, left, black] {$\bbover X$} ++ (0, -0.75);
\draw ([xshift=-0.25cm,yshift=-1.25cm]e) node[draw=black, fill=white, very thick, rounded corners, baseline=center] (f){$\overline{\eta}$};
\draw ([xshift=-0.8cm, yshift=-4.15cm]a) node[draw=black, fill=white, very thick, rounded corners, baseline=center] (g){$\overline{\xi}$};
\draw (g.south) --node[pos=0.5, left] {$\bbover b$}++ (0, -0.75);
\draw (f.south) --node[pos=0.5, left] {$\bbover a$}++ (0, -0.85);
\draw ([xshift=-0.25cm,yshift=-1.25cm]f) node[draw=black, fill=white, very thick, rounded corners, baseline=center] (h){$\nu^{-1}$};
\draw ([xshift=-0.05cm]h.south) --node[pos=0.5, left] {$\overline{ \bb a \otimes \bb b}$} ++ (0,-0.75);
\draw ([xshift=0.05cm]h.south) --++ (0,-0.75);
\draw ([yshift=-1.25cm]h) node[draw=black, fill=white, very thick, rounded corners, baseline=center] (i){$\overline{\alpha}$};
\draw (i.south) --node[pos=1, below] {$\bbover c$}++ (0, -0.5);
\draw [red, very thick] ([xshift=0.25cm]e.south) --node[pos=0.4, right, black] {$\bb X$}++ (0, -0.75);
\draw ([xshift=0.25cm,yshift=-1.25cm]e) node[draw=black, fill=white, very thick, rounded corners, baseline=center] (j){$\eta$};
\draw (j.south) --node[pos=0.7, left] {$\bb a$} ++ (0, -2.25);
\draw [red, very thick] (b) --++ (0, -3.25);
\draw ([yshift=-3.3cm]b) node[draw=black, fill=white, very thick, rounded corners, baseline=center] (k){$\xi$};
\draw (k.south) --node[pos=0.63, right] {$\bb b$}++ (0, -2.25);
\draw ([xshift=0.3cm,yshift=-2.5cm]j) node[draw=black, fill=white, very thick, rounded corners, baseline=center, inner xsep=0.3cm] (l){$\alpha$};
\draw (l.south) --node[pos=1, below] {$\bb c$} ++ (0, -0.5);
\end{tikzpicture}
\end{center}

By \cite[Example~11]{JP17}, $\psi_{\overline{c} \otimes c}$ is positive for $c \in \cC$.
Now we consider $\psi$ as a $\dagger$-functor of cyclic $\cC$ module $\rmC^*$ categories.
Then $\psi(\langle \eta | \eta \rangle_{\bb B}) = \overline{x}x$ for some $x \in \underline \End^{C^*}_{\cC}(\bb Y_{\bb C})(\bbover a \otimes \bb a)$.
Thus the above diagram is positive.

We apply Construction \ref{const:completion} to obtain a right $\bb C$ Hilbert $\rmC^*$ module object, denoted $\bb X \otimes_{\bb B} \bb Y$.
Note that elements of the form
\[
\begin{tikzpicture}[baseline={(current bounding box.center)}, thick]
\filldraw [black] (0,0) circle (0pt) coordinate (a);
\draw (a) to [out=180, in=-90] node[pos=0.5, left] {$\bb a$} ([xshift=-0.5cm, yshift=0.75cm]a.north) coordinate (b);
\draw (a) to [out=0, in=-90] node[pos=0.5, right] {$\bb b$} ([xshift=0.5cm, yshift=0.75cm]a.north) coordinate (c);
\node[above] at (a) {$\alpha$};
\draw (a) --node[pos=1, below] {$\bb c$}++ (0, -0.25);
\draw (b) node[draw=black, fill=white, very thick, rounded corners, baseline=center] (d){$\eta$};
\draw (c) node[draw=black, fill=white, very thick, rounded corners, baseline=center] (e){$\xi$};
\draw [red, very thick] (d.north) --node[pos=1, above, black] {$\bb X$}++ (0, 0.5) coordinate (f);
\draw [red, very thick] (e.north) to node[pos=1, above, black] {$\bb Y$} ([xshift=1cm]f);
\end{tikzpicture}
\quad \quad \text{ for all } a, b \in \Irr(\cC), \eta \in \bb X(a), \xi \in \bb Y(b), \alpha \in \cC(c, a \otimes b)
\]
are dense in $\bb X \otimes_{\bb B} \bb Y(c)$.
(This holds by a modification of \cite[Lemma~2]{JP17}.)
Thus we make $\bb X \otimes_{\bb B} \bb Y$ into an $\bb A-\bb C$ bimodule object by the composition
$$\bb A \to \underline \End^{C^*}_{\cC}(\bb X_{\bb B}) \to \underline \End^{C^*}_{\cC}(\bb X \otimes \bb Y_{\bb C}) \to \underline \End^{C^*}_{\cC}(\bb X \otimes_{\bb B} \bb Y_{\bb C})$$
where the first map comes from the data of $\bb X_{\bb B}$, the second map is given by inclusion, and the third map comes from Construction \ref{const:completion} applied to $\bb X \otimes \bb Y$.
These morphisms assemble into a natural transformation.
\end{const}

Now we are able to prove the main theorem of this section.

\begin{thm}
$\rmC^*$ algebra objects, Hilbert $\rmC^*$ bimodule objects, and intertwiners form a $\rmC^*$ 2-category $\CAlg_{\cC}$.
\end{thm}

\begin{proof}
To see that these form a 2-category, we construct associators and unitors, and show they satisfy the pentagon and triangle axioms.
For $\bb A-\bb B$, $\bb B-\bb C$, and $\bb C-\bb D$ Hilbert $\rmC^*$ bimodule objects $\bb X$, $\bb Y$, and $\bb Z$, the associator $\alpha_{\bb X, \bb Y, \bb Z}$ is given by applying the completion functor $\widetilde{\cdot}$ to the associator $\alpha_{\bb X, \bb Y, \bb Z}^{\Vect(\cC)}$.
This satisfies the pentagon axiom.
Its adjoint is its inverse, as can be seen from a similar computation to the one in the proof of Proposition \ref{prop:modulecat}.
Namely, that computation holds on a component-wise dense subset of an object in $\CAlg_{\cC}$, hence extends to the associator.
The unitors come from the fact that the relative product is balanced on components by construction, i.e. the expression $\eta \otimes \beta \otimes \lambda_{a, b}(\xi \otimes \gamma \otimes \zeta) - \rho_{a, b}(\eta \otimes \beta \otimes \psi_b(\xi)) \otimes \gamma \otimes \zeta$ for $\xi \in \bb X(a)$ and $\zeta \in \bb A(b)$ lies in $N^c$ (see Appendix \ref{app:completionconstruction}) for appropriately chosen $c$.
These are adjointable by construction of the relative product and since the unitors in $\CMod_{\bb A}$ are adjointable.

Now we show that these form a $\rmC^*$ 2-category.
The dagger structure on 2-morphisms is taking adjoints of adjointable operators.
By Proposition \ref{prop:modulecat} and our earlier remarks, the associators and unitors are unitary.
As the 1-categories of morphisms between objects have direct sums by Lemma \ref{lem:directsums}, endomorphism algebras in those categories are $\rmC^*$ algebras by Proposition \ref{prop:adjopcstar}.
\end{proof}

\section{Module Objects for $\rmW^*$ Algebra Objects} \label{section:WModObjs}
In this section we define $\WAlg_{\cC}$, the 2-category of $\rmW^*$ algebra objects, Hilbert $\rmW^*$ or Hilbert space bimodule objects, and normal adjointable intertwiners.
Theorem \ref{thm:B}, that Hilbert $\rmW^*$ and Hilbert space modules form equivalent categories, is a necessary ingredient for this definition.

\subsection{Hilbert $\rmW^*$ Module Objects}
\begin{defn} \label{defn:WModObjs}
Let $\bb M$ be a $\rmW^*$ algebra object in $\Vect(\cC)$.
We say $\bb X$ is a right $\bb M$ Hilbert $\rmW^*$ module object if $\bb X \in \CMod_{\bb M}$, there exists a predual $\bb X_*$, i.e. for all $c \in \cC$, a Banach space $\bb X_*(c)$ such that $(\bb X_*(c))^*$ is isometrically isomorphic to $\bb X(c)$, and the inner product is separately weak$^*$ continuous with respect to $\bb X_*$, i.e. for all $c, d \in \cC$ the following map $\bb X(d) \to \bb M(\overline d \otimes c)$ is weak$^*$ continuous for all $\eta \in \bf X(c)$:
\begin{center}
\begin{tikzpicture}[baseline={(current bounding box.center)}, thick]
\draw (0,0) node[draw=black, fill=white, very thick, rounded corners, baseline=center, inner xsep=0.3cm] (a){$\langle \cdot | \cdot \rangle_{\bb M}$};
\draw (a.north) --node[pos=1, above] {$\bb M$}++ (0, 0.5);
\draw [color=red] ([xshift=-0.35cm]a.south) --node[pos=0.4, left, color=black] {$\bbover X$}++ (0, -0.75);
\draw [color=red] ([xshift=0.35cm]a.south) --node[pos=0.4, right, color=black] {$\bb X$}++ (0, -0.75);
\draw ([xshift=-0.35cm, yshift=-1cm]a.south) node[draw=black, fill=white, very thick, rounded corners, baseline=center] (b){$\overline{\xi}$};
\draw ([xshift=0.35cm, yshift=-1cm]a.south) node[draw=black, fill=white, very thick, rounded corners, baseline=center] (c){$\eta$};
\draw (c.south) --node[pos=1, below] {$\bb c$}++ (0, -0.5) coordinate (d);
\draw (b.south) to node[pos=1, below] {$\bbover d$} ([xshift=-0.7cm]d);

\node[] (e) at ([xshift=-1.5cm, yshift=-0.5cm]a) {$\mapsto$};
\draw ([xshift=-1cm]e) node[draw=black, fill=white, very thick, rounded corners, baseline=center] (a){$\overline \xi$};
\draw (a.north)--node[pos=1, above] {$\bf X$}++ (0, 0.75);
\draw (a.south)--node[pos=1, below] {$\bf d$}++ (0, -0.75);
\end{tikzpicture}
\end{center}
We can replace this condition with a self-duality condition, a natural isomorphism
$$\bbover X \cong \left\{ \text{uniformly bounded maps } \bb X_{\bb M} \to \bb M_{\bb M} \right \}$$
with isometric components given by $\overline \eta \mapsto \langle \eta|$.
These conditions are equivalent by \cite[Lemma~8.5.4]{BLM04}, combined with the observation that for an $\bb M$ module $\bb X$, the action of $\bb M$ is captured by the actions of $\bb M(\overline c \otimes c)$ for all $c \in \cC$.
\end{defn}

\begin{rmk}
Also by \cite[Lemma 8.5.4]{BLM04}, \emph{the} predual $\bb X_*$ is unique up to unique isometric isomorphism.
Suppose that $\bb Y$ is also a predual for $\bb X$ in the sense of Definition \ref{defn:WModObjs}.
The lemma gives us that $\bb X_*(c)$ is isometrically isomorphic to $\bf Y(c)$.
Semisimplicity of $\cC$ enforces that this isomorphism is natural.
\end{rmk}

To make sense of the separately weak$^*$ continuous condition of Definition \ref{defn:WModObjs}, $\bb M(c)$ must have a predual, which it is does as $\bb M$ is a $\rmW^*$ algebra object.
More generally, we have the following example.

\begin{ex}
A $\rmW^*$ algebra object $\bb M$ is a right Hilbert $\rmW^*$ module over itself.
We saw in Example \ref{ex:overself} that $\bb M \in \CMod_{\bb M}$.
As $\bb M$ is a $\rmW^*$ module object, $\cM_{\bb M}$ is a $\rmW^*$ category.
Thus $\bb M(c) \cong \cM_{\bb M}(c \otimes \bb M, \bb M)$ has a predual.
Combining this with the fact that conditional expectations are normal implies the inner product on $\bb M$ is separately weak$^*$ continuous.
\end{ex}

As in the $\rmC^*$ case, we also have this example.

\begin{ex} \label{ex:tensorbycompacts}
For $\bb X \in \WMod_{\bb M}$, $\bb c \otimes \bb X \in \WMod_{\bb M}$ as well.
By Example \ref{ex:cyclicmodulecat}, $\bb c \otimes \bb X \in \CMod_{\bb M}$.
To provide a predual $(\bb c \otimes \bb X)_*$, it suffices to produce components over simple objects of $\cC$.
For $a \in \Irr(\cC)$, $(\bb c \otimes \bb X)(a) = \bigoplus \limits_{d \in \Irr(\cC)} \cC(d, c \otimes a) \otimes \bb X(a)$.
So $(\bb c \otimes \bb X)_*(a) = \bigoplus \limits_{d \in \Irr(\cC)} \cC(d, c \otimes a) \otimes \bb X_*(a)$ is a predual for $(\bb c \otimes \bb X)(a)$ as $\cC(d, c \otimes a)$ is a finite dimensional Hilbert space and $c \otimes a$ can be written as a finite direct sum of simple objects in $\cC$.
The inner product on $\bb c \otimes \bb X$ is separately weak$^*$ continuous as $\cM_{\bb M}$ is a $\rmW^*$ category.
\end{ex}

We move towards showing $\WMod_{\bb M}$ is a $\rmW^*$ module category.
We must first define morphisms.

\begin{defn}
For $\bb X, \bb Y \in \WMod_{\bb M}$, $\alpha \in \Vect(\cC)(\bb X, \bb Y)$ is a morphism in $\WMod_{\bb M}$ if $\alpha$ is a morphism in $\CMod_{\bb M}$.
\end{defn}

Note that this implies $\alpha_c$ is weak$^*$ continuous for all $c \in \cC$ by \cite[Corollary~8.5.8]{BLM04}, upon viewing $\bb X(c)$ and $\bf Y(c)$ as $\bb M(1_{\cC})$ Hilbert $\rmW^*$ modules.

\begin{prop} \label{prop:WModWStar}
The category $\WMod_{\bb M}$ is a $\cC$ module $\rmW^*$ category.
\end{prop}

\begin{proof}
That $\WMod_{\bb M}$ is a $\cC$ module category is similar to the proof of Proposition \ref{prop:modulecat}.
Normality of the associators and unitors holds easily.
To see $\WMod_{\bb M}$ admits direct sums, we proceed as follows.
By Lemma \ref{lem:directsums}, $\CMod_{\bb M}$ admits direct sums.
The direct sum $\bb X \oplus \bb Y$ has predual $\bb X_* \oplus \bb Y_*$, and the $\bb M$-valued inner product it carries is separately weak$^*$ continuous.
By the form of the predual of the direct sum, the morphisms in the direct sum diagram have weak$^*$ continuous components.
One can quickly check the universal property from there.

Finally, we check that endomorphism algebras in $\WMod_{\bb M}$ are von Neumann algebras.
For $\bb X \in \WMod_{\bb M}$, we pick a convenient faithful $*$-representation of $N = \WMod_{\bb M}(\bb X, \bb X)$.
For $c \in \Irr(\cC)$, let $H(c)$ be the completion of $\bb X(c)$ with respect to the inner product $\langle \eta| \xi \rangle$ given by $\tau$ applied to the diagram
\begin{center}
\begin{tikzpicture}[baseline={(current bounding box.center)},thick]
\draw (0, 0) node[draw=black,fill=white,very thick,rounded corners,baseline=center] (a){$\langle \cdot | \cdot \rangle_{\bb M}$};
\draw [red, very thick] ([xshift=-0.3cm, yshift=-0.75cm]a.south) --node[pos=0.5, left, black] {$\bbover X$}++(0, 0.75);
\draw [red, very thick] ([xshift=0.3cm, yshift=-0.75cm]a.south) --node[pos=0.5, right, black] {$\bb X$}++(0, 0.75);
\draw (a.north) --node[pos=1, above, black] {$\bb M$}++(0, 0.5);

\draw ([xshift=-0.3cm, yshift=-1.35cm]a) node[draw=black,fill=white,very thick,rounded corners,baseline=center] (b){$\overline \eta$};

\draw ([xshift=0.3cm, yshift=-1.35cm]a) node[draw=black,fill=white,very thick,rounded corners,baseline=center] (c){$\xi$};
\filldraw ([xshift=0.3cm, yshift=-0.3cm]b.south) circle (0pt) coordinate (d);
\draw (b.south) to [out=-90, in=180] node[pos=0.5, left] {$\bbover c$} (d);
\draw (c.south) to [out=-90, in=0] node[pos=0.5, right] {$\bb c$} (d);
\end{tikzpicture}
\end{center}
Let $H = \bigoplus \limits_{c \in \Irr{\cC}} H(c)$.
We define the representation by $\beta \mapsto \bigoplus \beta_c$, where $\beta_c$ acts on $H(c)$ by $\eta \mapsto \beta_c(\eta)$ for $\eta \in \bfX(c)$, and extending by continuity.
As $\|\beta\| < \infty$, the image of $\beta$ extends to an operator on $H$.

Now suppose the net $(\alpha^i) \subset N$ converges strongly in $B(H)$ to $\alpha$.
We want to show that $\alpha \in N$.
First, we claim that $\bb X(f) \alpha_c = \alpha_d \bb X(f)$ for $f \in \cC(d, c)$, where by $\alpha_c$, we mean the SOT limit of $\alpha^i_c$.
For $\eta \in \bb X(c)$, observe the computation
$$\bb X(f)\alpha_c(\eta) = \bb X(f) (\lim \alpha^i_c(\eta)) \stackrel{(*)}{=} \lim (\bb X(f) \alpha^i_c(\eta)) \stackrel{(**)}{=} \lim (\alpha^i_d \bb X(f)(\eta)) = \alpha_d \bb X(f)(\eta)$$
where $(*)$ holds as $\bb X$ is a $\rmW^*$ algebra object and $(**)$ holds by naturality of $\alpha^i$.
Since $\bb X(c)$ is dense in $H(c)$, $\alpha$ is a natural transformation.
Observe $\alpha$ has an adjoint ($\alpha_c^* = \lim (\alpha^i_c)^*$ suffices).
To see that $\alpha$ is an adjointable operator, we show that $\|\alpha\| < \infty$.
$$\sup \limits_{c \in \Irr(\cC)} \|\alpha_c\| = \sup \limits_{c \in \Irr(\cC)} \| \lim \limits_i \alpha_c^i\| \stackrel{(*)}{\leq} \sup \limits_{c \in \Irr(\cC)} \lim \limits_i \| \alpha_c^i\| \leq \lim \limits_i \|\alpha^i\| \stackrel{(**)}{<} \infty$$
where $(*)$ and $(**)$ hold as $(\alpha^i)$ is weak$^*$ continuous, hence norm-bounded.
\end{proof}

Thus $\underline \End^{\rmW^*}_{\cC}(\bb X_{\bb M})$, defined analogously to Definition \ref{defn:intend}, is a $\rmW^*$ algebra object.

\subsection{Hilbert Space Module Objects}
\begin{defn}
For $\bb H \in \Hilb(\cC)$, we say $\bb H$ is a right Hilbert space module object for the $\rmW^*$ algebra object $\bb M$ if there exists a normal $*$-algebra natural transformation $\bb M^{\op} \to \bb B(\bb H)$, or equivalently $\bb M \to \bb B (\bb H)^{\op}$.
\end{defn}
~
\begin{ex}~
\begin{itemize}
\item $\bb H$ over $\bb B (\bb H)^{\op}$.
\item $L^2 \bb M$ over $\bb M$ (by \cite[Corollary~10]{JP17}).
\item $\bb c \otimes \bb H$ over $\bb M$ where $\bb H$ is a right Hilbert space module over $\bb M$.
\end{itemize}
\end{ex}

\begin{defn}
For $\bb H, \bb K \in \HilbMod_{\bb M}$, we say that $\alpha \in \HilbMod_{\bb M}(\bb H, \bb K)$ if $\alpha$ is a uniformly bounded right $\bb M$ module morphism.
\end{defn}

\begin{prop}
The category $\HilbMod_{\bb M}$ is a $\cC$ module $\rmW^*$ category.
\end{prop}

\begin{proof}
The module action is given by the tensor product in $\Hilb(\cC)$.
So the associator and unitors are given by the appropriate morphisms in $\Hilb(\cC)$.
(These are normal as $\Hilb(\cC)$ is a $\rmW^*$ category.)
It is clear that the category admits direct sums as $\Hilb(\cC)$ does.
To finish the proof, it suffices to check that $X = \HilbMod_{\bb M}(\bb H, \bb H)$ is a von Neumann algebra.
First, $X$ is a $*$-algebra under the adjoint of morphisms in $\Hilb(\cC)$.
That $X$ is SOT closed follows from the fact that the SOT limit of module maps is a module map, and the argument in Proposition \ref{prop:WModWStar} for $\|\alpha\| < \infty$.
\end{proof}

\subsection{The Equivalence of Hilbert Space and Hilbert $\rmW^*$ Module Objects}
In this section we show that the 1-categories exposited in the previous subsections are equivalent.
We begin by defining bimodules for both Hilbert $\rmW^*$ module objects and Hilbert space module objects.

\begin{defn}
For $\bb X \in \WMod_{\bb M}$, $\bb X$ is a $\bb N-\bb M$ Hilbert $\rmW^*$ bimodule object if there exists a unital normal $*$-algebra natural transformation $\phi \colon \bb N \to \underline \End^{\rmW^*}_{\cC}(\bb X_{\bb M})$.
\end{defn}

\begin{ex}
A $\rmW^*$ algebra object $\bb M$ is a Hilbert $\rmW^*$ bimodule object over itself.
As in Example \ref{ex:algebraoverself}, this is an equivalence of $\cM_{\bfM}$ and $\cM_{\bfM}^{\rmW^*}$.
\end{ex}

As with Hilbert $\rmC^*$ bimodule objects, there is a relative tensor product for Hilbert $\rmW^*$ bimodule objects.
Let ${}_{\bfL} \bfX_{\bfM}$ and  ${}_{\bfM} \bfY_{\bfN}$ be Hilbert $\rmW^*$ bimodule objects and consider $\bfX \otimes_{\bfM} \bfY_{\bfN}$ as in the $\rmC^*$ case.
This has a predual by \cite[8.5.40]{BLM04}.
Moreover, the inner product on $\bf X \otimes_{\bfM} \bfY_{\bfN}$ is separately weak$^*$ continuous, and the $*$-algebra natural transformation $\bb L \to \underline \End^{\rmW^*}_{\cC}(\bb X \otimes_{\bb M} \bb Y_{\bb N})$ is normal, as $\WMod_{\bfM}$ and $\WMod_{\bfN}$ are $\rmW^*$ category.

On the other hand, we define Hilbert space bimodule objects as follows.
\begin{defn}
For $\bb H$ a left $\bb M$ Hilbert space module object and a right $\bb N$ Hilbert space module object, we say that $\bb H$ is an $\bb M - \bb N$ Hilbert space bimodule object if the actions of $\bb M$ and $\bb N$ on $\bb H$ commute.\footnote{One can interpret this definition in terms of commutants using the language of \cite[Definition~44]{JP17}.}
\end{defn}

To induce a monoidal structure on $\HilbMod_{\bb M}$, we pass through the monoidal structure on Hilbert $\rmW^*$ bimodule objects.
The following functors allows us to achieve this goal.

\begin{const}
We construct the functor $- \otimes_{\bb M} L^2 \bb M$.
For $\bb X \in \WMod_{\bb M}$, we define $\bb X \otimes_{\bb M} L^2 \bb M$ componentwise as the relative product (see Construction \ref{const:relprod}) with respect to the following densely-defined inner product.
For $\eta_i \otimes \alpha_i \otimes \xi_i \in \bb X \otimes \widehat{\bb M}$, let $\langle \eta_1 \otimes \alpha_1 \otimes \xi_1 | \eta_2 \otimes \alpha_2 \otimes \xi_2 \rangle$ be $\tau$ applied to the diagram:
\begin{center}
\begin{tikzpicture}[baseline={(current bounding box.center)}, thick]
\draw (0,0) node[draw=black, fill=white, very thick, rounded corners, baseline=center, inner xsep=0.3cm] (a){$\langle \cdot | \cdot \rangle_{\bb M}$};
\draw [color=red] ([xshift=-0.35cm]a.south) --node[pos=0.4, left, color=black] {$\bbover X$}++ (0, -0.75);
\draw [color=red] ([xshift=0.35cm]a.south) --node[pos=0.4, right, color=black] {$\bb X$}++ (0, -0.75);
\draw ([xshift=-0.35cm, yshift=-1cm]a.south) node[draw=black, fill=white, very thick, rounded corners, baseline=center] (b){$\overline{\eta_1}$};
\draw ([xshift=0.35cm, yshift=-1cm]a.south) node[draw=black, fill=white, very thick, rounded corners, baseline=center] (c){$\eta_2$};
\filldraw [black] ([xshift=-0.5cm, yshift=-0.75cm]b.south) circle (0pt) coordinate (d);
\node[above] at (d) {$\overline{\alpha_1}$};
\filldraw [black] ([xshift=0.5cm, yshift=-0.75cm]c.south) circle (0pt) coordinate (e);
\node[above] at (e) {$\alpha_2$};
\draw ([xshift=-1cm, yshift=0.05cm]b) node[draw=black, fill=white, very thick, rounded corners, baseline=center] (f){$\overline \xi_1$};
\draw ([xshift=1cm, yshift=0.05cm]c) node[draw=black, fill=white, very thick, rounded corners, baseline=center] (g){$\xi_2$};
\draw (d) to [out=180, in=-90] (f.south);
\draw (d) to [out=0, in=-90] (b.south);
\draw (e) to [out=0, in=-90] (g.south);
\draw (e) to [out=180, in=-90] (c.south);
\filldraw ([xshift=0.35cm, yshift=-1.5cm]b.south) circle (0pt) coordinate (h);
\filldraw ([yshift=1cm]a.north) circle (0pt) coordinate (i);
\draw (h) to [out=180, in=-90] (d);
\draw (h) to [out=0, in=-90] (e);
\draw (f.north) --node[pos=0.5, left] {$\widehat{\bb M}$}++ (0,1);
\draw (g.north) --node[pos=0.5, right] {$\widehat{\bb M}$}++ (0,1);
\draw ([yshift=1cm]f.north) to [out=90, in=180] (i);
\draw ([yshift=1cm]g.north) to [out=90, in=0] coordinate (j) (i);
\draw (a.north) to [out=90, in=180] node[pos=0.3, left] {$\bb M$} (j);
\draw (i) --node[pos=1, above] {$\widehat{\bb M}$}++ (0, 0.5);
\end{tikzpicture}
\end{center}
By construction, this is an object in $\Hilb(\cC)$.
Since $\bb M$ is tracial, we have a normal right action of $\bb M$ on $L^2 \bb M$, hence on $\bb X \otimes_{\bb M} L^2 \bb M$.
For $\alpha \in \WMod_{\bb M}(\bb X, \bb Y)$, we obtain a densely-defined morphism $\alpha \otimes \id \colon \bb X \otimes L^2\bb M \to \bb Y \otimes_{\bfM} L^2\bb M$, which extends by continuity to a morphism on $\bfX \otimes_{\bfM} L^2\bfM$.
By construction, $\alpha \otimes \id$ is a morphism of Hilbert space module objects.
From this, one can quickly check that  $- \otimes_{\bb M} L^2\bb M$ is a normal dagger functor.
\end{const}

On the other hand, we can go from $\HilbMod_{\bb M}$ to $\WMod_{\bb M}$ by taking (left) bounded vectors.

\begin{const}
Here we construct the functor $D(-)$.
For $\bb H \in \HilbMod_{\bb M}$, we define $D(\bb H)$ to be $\HilbMod_{\bb M}(L^2\bb M, \bf H)$.
The right action of $\bb M$ on $D(\bb H)$ is given by precomposition with the left action of $\bb M$ on $L^2\bb M$.
To make $D(\bb H)$ into a Hilbert $\rmW^*$ module object, we equip it with an $\bfM$-valued inner product.
For $\alpha, \beta \in D(\bb H)$, we naturally identify $\alpha^\dagger_c \beta_c$ with an element of $\bb M(c)$ to provide this structure by bootstrapping the equivalence of Hilbert $\rm W^*$ modules and Hilbert space modules for \emph{honest} von Neumann algebras.

First consider $f \in \HilbMod_{\bb M} (L^2 \bb M, L^2 \bb M)$.
We can describe a bounded map of Hilbert spaces $f_c \colon L^2 \bb M(c) \to L^2 \bb M(c)$ by restricting the domain to a dense set, say $\bb M(c)$.
Now, $\bb M(c)$ is a right $\bb M(1_{\cC})$ Hilbert $\rmW^*$ module, so $D(\bb M(c) \otimes_{\bb M(1_{\cC})} L^2(\bb M(1_{\cC}))$ is isometrically isomorphic to $\bb M(c)$.
Since bounded $\bb M(1_{\cC})$ maps, such as $f_c$, take bounded vectors to bounded vectors, we observe that $f_c(\bb M(c)) \subset \bb M(c)$.
This lets us write $f$ as a right $\bb M$ linear endomorphism of $\bb M$.
Clearly $f$ is normal.
So we can represent $f_c$ as some element of $\bb M$.
Where $f = \alpha^\dagger \beta$ for $\alpha, \beta \in D(\bb H)$, this provides us with the desired identification.
Note that \ref{H2} and \ref{H4} follow from the fact that $\WMod_{\bb M}$ is a $\rmW^*$ category.

This $\bb M$-valued inner product makes $D(\bb H)$ into a Hilbert $\rmW^*$ module.
We see this by showing $D(\bb H)$ is self-dual.
The map $\overline{\alpha_c} \mapsto \langle \alpha_c|$ is injective as $\langle \cdot | \cdot \rangle_{\bb M}$ is non-degenerate.
Take a projective basis $\{e_i^c\}$ for $D(\bb H)(c)$.
We can describe $\beta_c \in D(\bb H)(c)$ by the image of the projective basis under $\beta_c$.
Say $\beta_c(e_i^c) = \eta_i$.
Then there exists $\alpha_c$ such that $\langle \alpha_c | e_c^i \rangle = \eta_i$ as $D(\bfH)(c)$ is an $\bfM(1_{\cC})$ Hilbert $\rmW^*$ module.
Since we can perform this procedure for all $c \in \cC$, and $\beta$ is a natural transformation, $\alpha_c$ is a component of the natural transformation $\alpha$.
Note $\langle \alpha_c | = \beta_c$.
Moreover, $\|\alpha_c\| = \|\langle \alpha_c| \|$; the left inequality follows trivially, and the right one follows from Cauchy-Schwarz.
Therefore $D(\bb H)$ is self-dual.

We make $D(-)$ into a functor by taking $D(\alpha)$ to be postcomposition for $\alpha \in \HilbMod_{\bb M}(\bb H, \bb K)$.
The adjoint of $D(\alpha)$ is the same as its adjoint in $\Hilb(\cC)$, and $D(\alpha)_c$ is weak$^*$ continuous for all $c \in \cC$ as $\Hilb(\cC)$ is a $\rmW^*$ category.
Thus $D(-)$ is a normal dagger functor.
\end{const}

Observe that if $\bb H$ is an $\bb N-\bb M$ Hilbert space bimodule object, then $D(\bb H)$ has a left $\bb N$ action given by postcomposition with the left $\bb N$ action on $\bb H$.
Since this commutes with the right action of $\bb M$, this describes a normal $*$-algebra natural transformation $\bb N \to \underline \End^{\rmW^*}_{\cC}(D(\bb H)_{\bb M})$.

Now we prove Theorem \ref{thm:B}.

\begin{thm}
The functors $$- \otimes_{\bb M} L^2 \bb M \colon \WMod_{\bb M} \to \HilbMod_{\bb M}$$ and $$D(-) \coloneq \HilbMod_{\bb M}(L^2 \bb M, -) \colon \HilbMod_{\bb M} \to \WMod_{\bb M}$$ witness a $\dagger$-equivalence of $\rmW^*$ categories.
\end{thm}

\begin{proof}
We show the following isomorphisms in $\HilbMod_{\bb M}$ and $\WMod_{\bb M}$, respectively:
$$ D(\bb H) \otimes_{\bb M} L^2 \bb M \cong \bb H \quad\quad\quad D(\bb X \otimes_{\bb M} L^2 \bb M) \cong \bb X$$

For the first, observe that the set $D(\bb H)(c) \otimes \Omega$ is dense in $D(\bb H) \otimes_{\bb M} L^2\bb M(c)$.
So it suffices to define the first map for $\alpha_c \otimes \Omega$, where $\alpha \in D(\bb H)$.
We send $\alpha_c \otimes \Omega$ to $\alpha_c(x \Omega)$ for $x = \sum e_i^c$, and $\{e_i^c\}$ a projective basis for $\bb M(c)$.
This is clearly right $\bb M$-linear and natural.
It is isometric as 

\begin{center}
\begin{tikzpicture}[baseline={(current bounding box.center)}, thick]
\node[] (a) at (0,0) {$\|\langle \alpha_c | \alpha_c \rangle_{\bb M} \| = \langle x \Omega | x \Omega \rangle_{L^2\bb M(c)} || \langle \alpha_c | \alpha_c \rangle_{\bfM}||^2 = $};

\draw ([xshift=5.25cm, yshift=-1.25cm]a) node[draw=black, fill=white, very thick, rounded corners, baseline=center] (a){$\alpha$};
\draw ([yshift=1.25cm]a) node[draw=black, fill=white, very thick, rounded corners, baseline=center] (b){$\alpha^*$};
\draw ([yshift=1.25cm]b) node[draw=black, fill=white, very thick, rounded corners, baseline=center] (c){$x^* \Omega$};
\draw ([yshift=1.25cm]c) node[draw=black, fill=white, very thick, rounded corners, baseline=center] (d){$\varphi_c$};
\draw ([yshift=-1.25cm]a) node[draw=black, fill=white, very thick, rounded corners, baseline=center] (e){$x \Omega$};
\draw ([xshift=-1cm]d.north) coordinate (f) to node[pos=0.5, left] {$\bbover c$} ([xshift=-1cm] e.south)coordinate (g);
\filldraw ([xshift=0.5cm, yshift=0.5cm]f) circle (0) coordinate (h);
\filldraw ([xshift=0.5cm, yshift=-0.5cm]g) circle (0) coordinate (i);
\draw (d.north) to [out=90, in=0] (h);
\draw (h) to [out=180, in=90] (f);
\draw (e.south) to [out=-90, in=0] (i);
\draw (i) to [out=180, in=-90] (g);
\draw (a.north) to node[pos=0.5, right] {$\bb H$} (b.south);
\draw (b.north) to node[pos=0.5, right] {$L^2 \bb M$} (c.south);
\draw (c.north) to node[pos=0.5, right] {$\bb c$} (d.south);
\draw (e.north) to node[pos=0.5, right] {$L^2 \bb M$} (a.south);
\end{tikzpicture}
\end{center}

The first equality holds by traciality.
Moreover, this map has dense range, since $x$ is the sum of a projective basis.
So we have our first unitary isomorphism.

For the second, it suffices to consider the component of the map over a simple object $c \in \cC$.
For $\eta \in \bb X(c)$, we define a natural transformation $\bb X \to D(\bb X \otimes_{\bb M} L^2 \bb M)$ by 
\begin{center}
\begin{tikzpicture}[baseline={(current bounding box.center)}, thick]
\draw (0,0) node[draw=black, fill=white, very thick, rounded corners, baseline=center] (a){$x \Omega$};
\draw (a.south) --node[pos=1, below] {$\bb c$}++ (0, -0.5);
\draw (a.north) --node[pos=1, above] {$\widehat{\bb M}$}++ (0, 0.5);

\node[] (b) at ([xshift=1cm]a) {$\mapsto$};

\draw ([xshift=0.75cm]b) node[draw=black, fill=white, very thick, rounded corners, baseline=center] (a){$\eta$};
\draw (a.north) --node[pos=1, above] {$\bb X$}++ (0, 0.5) coordinate (c);
\draw (a.south) --node[pos=1, below] {$\bb c$}++ (0, -0.5);
\filldraw ([xshift=0.5cm]a.east) circle (2pt) coordinate (b);
\draw (b) to node[pos=1, above] {$\widehat{\bb M}$} ([xshift=0.75cm]c);
\end{tikzpicture}
\end{center}
on the $c$ component and as 0 on all other components over simple objects.
This is natural as $c$ is simple and right $\bb M$ linear by the construction of the relative product.
Moreover, the assignment of $\eta$ to this natural transformation is itself manifestly natural.
The map is isometric $\bb X(c) \to D(\bb X \otimes_{\bb M} L^2 \bb M)(c)$ by faithfulness of the trace $\tau$ on $\bb M(1_{\cC})$.
We have dense range by the observation:
\begin{center}
\begin{tikzpicture}[baseline={(current bounding box.center)}, thick]
\draw (0,0) node[draw=black, fill=white, very thick, rounded corners, baseline=center] (a){$x \Omega$};
\draw (a.south) --node[pos=1, below] {$\bb c$}++ (0, -1) coordinate (y);
\draw (a.north) --node[pos=1, above] {$\widehat{ \bb M}$}++ (0, 0.5) coordinate (z);

\node[] (b) at ([xshift=1cm, yshift=-0.35cm]a) {$\mapsto$};

\draw ([xshift=1.25cm, yshift=0.35cm]b) node[draw=black, fill=white, very thick, rounded corners, baseline=center] (a){$\alpha_c(x \Omega)$};
\draw (a.south) to node[pos=1, below] {$\bb c$} ([xshift=2.25cm]y);
\draw (a.north) to node[pos=1, above] {$\bb X \otimes_{\bb M} L^2 \bb M$} ([xshift=2.25cm]z);

\node[] (z) at ([xshift=1.25cm, yshift=-0.35cm]a) {$=$};

\draw ([xshift=1.25cm, yshift=0.35cm]z) node[draw=black, fill=white, very thick, rounded corners, baseline=center] (a){$\alpha_{1_{\cC}}(\Omega)$};
\draw (a.north) --node[pos=1, above] {$\bb X$}++ (0, 0.5) coordinate (e);
\draw ([xshift=2cm, yshift=-0.75cm]a) node[draw=black, fill=white, very thick, rounded corners, baseline=center] (b){$j_c(x)$};
\draw(b.south) --node[pos=1, below] {$\bb c$}++ (0, -0.25);
\filldraw ([xshift=1cm]a.south) circle (2pt) coordinate (c);
\draw (c) to node[pos=1, above] {$\widehat{\bb M}$} coordinate (d) ([xshift=1cm]e);
\draw (b.north) to [out=90, in=0] node[pos=0.5, above] {$\bb M$} (d);
\end{tikzpicture}
\end{center}
for $\alpha \in D(\bb X \otimes_{\bb M} L^2 \bb M)$, so we have a unitary isomorphism.
\end{proof}

Using the monoidal structure on Hilbert $\rmW^*$ bimodules and this equivalence of categories, we induce a monoidal structure on Hilbert space bimodules by $$\bb{H} \boxtimes_{\bb M} \bb K \coloneq D(\bb H) \otimes_{\bb M} D(\bb K) \otimes_{\bb N} L^2 \bb N$$ completed component-wise with respect to the inner product $$\langle f_1 \otimes \alpha_1 \otimes g_1 \otimes \beta_1 \otimes \eta_1 | f_2 \otimes \alpha_2 \otimes g_2 \otimes \beta_2 \otimes \eta_2\rangle \coloneq \langle \langle g_1 | g_2 \rangle_{D(\bb K)} \langle f_1 | f_2 \rangle_{D(\bb H)} \eta_1 | \eta_2 \rangle_{L^2 \bb M(a)}$$ where $\eta_i \in L^2 \bb M(a)$.
This monoidal structure is a $\cC$-graded Connes' fusion tensor product.
Thus, we see that:

\begin{cor} \label{cor:WAlgwelldef}
The following 2-categories are $\dagger$-equivalent:
\[
  \left\{ \begin{array}{c}
    \rmW^* \text{ algebra objects} \\
    \text{Hilbert } \rmW^* \text{ bimodules} \\
    \text{adjointable intertwiners}
  \end{array}\right\}
  \text{ and }
  \left\{ \begin{array}{c}
    \rmW^* \text{ algebra objects} \\
    \text{Hilbert space bimodules} \\
    \text{adjointable intertwiners}
  \end{array}\right\}.
\]
Thus we can call either one $\WAlg_{\cC}$.
\end{cor}
As in the $\rmC^*$ case, we have that $\WAlg_{\cC}$ is a $\rmW^*$ 2-category.
Since the left and right models of $\WAlg_{\cC}$ are equivalent via normal dagger functors, it suffices to check this for the left model.
By Proposition \ref{prop:WModWStar}, it remains to see that the relative tensor product is separately weak$^*$ continuous.
This check occurs in the category $\WAlg$, or $\WAlg_{\Hilb_{\rm fd}}$, which is a $\rmW^*$ 2-category.
So $\WAlg_{\cC}$ is indeed a $\rmW^*$ 2-category.

\section{An Operator Algebraic Ostrik's Theorem}
In this section, we prove the main result of the paper,  Theorem \ref{thm:A}.

\begin{thm}
There exist $\dagger$-equivalences of 2-categories $\CAlg_{\cC} \cong \CMod_{\cC}$ and $\WAlg_{\cC} \cong \WMod_{\cC}$.
\end{thm}

\begin{proof}
We begin with the first equivalence.
The 2-functor $\phi^{\rmC^*}$ is defined on objects by $\bb A \mapsto \cM_{\bb A}$, on 1-morphisms by ${}_{\bb A} \bb X_{\bb B} \mapsto - \otimes_{\bb A} \bb X_{\bb B}$, and on 2-morphisms by sending $\alpha \colon {}_{\bb A} \bb X_{\bb B} \to {}_{\bb A} \bb Y_{\bb B}$ to $- \otimes \alpha$.
It is easy to see that this is a $\dagger$ 2-functor.
By Example \ref{ex:algebraoverself}, combined with noting that morphisms $\cM_{\bfA}(c \otimes \bfA, d \otimes \bf A)$ are equivalent to elements of $\bfA(\overline d \otimes c)$, we can exchange $\cM_{\bb A}$ for $\cM^{\rmC^*}_{\bb A}$.
Thus $\phi^{\rmC^*}$ is essentially surjective on objects by the definition of a $\rmC^*$ algebra object.

Next, we show that $\CAlg_{\cC}(\bb A, \bb B) \cong \CMod_{\cC}(\cM_{\bb A}^{\rmC^*}, \cM_{\bb B}^{\rmC^*})$ as $\rmC^*$ categories.
By Ostrik's theorem, for $F \in \CMod_{\cC}(\cM_{\bb A}^{\rmC^*}, \cM_{\bb B}^{\rmC^*})$, there exists $\bb X \in \Alg_{\cC}(\bb A, \bb B)$ such that $\gamma$ provides an isomorphism $F \cong - \otimes_{\bb A} \bb X_{\bb B}$ in $\Mod_{\cC}(\cM_{\bb A}^{\rmC^*}, \cM_{\bb B}^{\rmC^*})$.
As $F(\bb A) \cong \bb X_{\bb B}$, and the target of $F$ is $\cM_{\bb B}$, $X_{\bb B} \in \CMod_{\bb B}$.
We see $\bb X \in \CMod_{\cC}(\cM_{\bb A}^{\rmC^*}, \cM_{\bb B}^{\rmC^*})$ via the map on $\bb A$ induced by $F$; observe $F(\bb A) = \underline \End^{\rmC^*}_{\cC}(\bb X_{\bb B})(a)$.
Then $\gamma$ is \emph{actually} in $\CMod_{\cC}(\cM_{\bb A}^{\rmC^*}, \cM_{\bb B}^{\rmC^*})$ as $F$ is a module functor: 
\[\|\gamma\| = \sup \limits_{c \in \Irr(\cC)} \|\gamma_c\| = \|\gamma_{1_{\cC}}\|.\]
That is, the $c$ component of $\gamma$ is a morphism $\gamma_c \colon c \otimes F(\bf A) \to \bf c \otimes \bfX$ as $F(\bf c \otimes \bfA) \cong \bf c \otimes F(\bf A)$ and $(\bf c \otimes \bf A) \otimes_{\bf A} \bf X \cong \bf c \otimes \bf X$.
Thus $\phi^{\rmC^*}$ is essentially surjective on 1-morphisms.

It is clear that $\phi^{\rmC^*}$ is injective on 2-morphisms.
Surjective on 2-morphisms holds like so.
By Ostrik's theorem, for all $\alpha \colon - \otimes_{\bb A} \bb X_{\bb B} \to - \otimes_{\bb A} Y_{\bb B} \in \CMod_{\cC}(- \otimes_{\bb A} \bb X_{\bb B}, - \otimes_{\bb A} \bb Y_{\bb B})$, there exists $f \in \Alg_{\cC}({}_{\bb A} \bb X_{\bb B}, {}_{\bb A} \bb Y_{\bb B})$ such that $\alpha = - \otimes f \in \Mod_{\cC}(- \otimes_{\bb A} \bb X_{\bb B}, - \otimes_{\bb A} Y_{\bb B})$.
Then $f \in \CAlg({}_{\bb A} \bb X_{\bb B}, {}_{\bb A} \bb X_{\bb B})$ as $f = \alpha_{\bb A}$, so $\phi^{\rmC^*}$ is surjective on 2-morphisms.

The proof of the second equivalence is similar.
The only interesting detail is seeing that $\gamma \colon F \to - \otimes_{\bfA} \bfX_{\bfB}$ is an isomorphism in $\WMod_{\cC}(\cM_{\bfA}^{\rmW^*}, \cM_{\bfB}^{\rmW^*})$.
Our notation is the notation above, but with $\rmW^*$ algebra objects $\bfA$ and $\bfB$.
Since $\gamma \in \CMod_{\cC}(\cM_{\bfA}^{\rmW^*}, \cM_{\bfB}^{\rmW^*})$, this follows from \cite[Corollary~8.5.10]{BLM04}.
\end{proof}

\begin{cor}
The 2-categories $\CMod_{\cC}$ and $\WMod_{\cC}$ are $\rmC^*$ and $\rmW^*$ 2-categories, respectively.
\end{cor}

\begin{rmk} \label{rmk:HNPgeneralized}
On the surface, it may appear that this result does not quite generalize \cite[Theorem~7.4]{HNP24}.
Namely, as of yet we have not addressed the (Cauchy) completeness and bi-involutive conditions.
We do so now.
By definition, the objects of $\WMod_{\cC}$ are finitely generated complete $\rm W^*$ categories.
Moreover, if one uses the version of $\WAlg_{\cC}$ where 1-morphisms are Hilbert space bimodule objects, then Theorem \ref{thm:A} is a bi-involutive equivalence.
However, note that Theorem \ref{thm:B}, and hence Corollary \ref{cor:WAlgwelldef}, is \emph{not} a bi-involutive equivalence.
We leave the verification of these statements to the interested reader.
\end{rmk}

\appendix

\section{The Completion Construction for Pre-Hilbert $\rmC^*$ Module Objects} \label{app:completionconstruction}

\begin{const} \label{const:completion} 
The following procedure builds a Hilbert $\rmC^*$ module object $\widetilde{\bb X}$ out of a right $\bb A$ module object $\bb X$ with a right $\bb A$-valued pre-inner product.
We denote the action of $\bb A$ on $\bb X$ by $\lambda$ and the pre-inner product by $\langle \cdot | \cdot \rangle_{\bb A}$.
To describe $\lambda$, the data of $\lambda_{a, b} \colon \bb X(a) \otimes \bb A(b) \to \bb X(a \otimes b)$ for $a, b \in \Irr(\cC)$  suffices.
Let $$N^c \coloneq \{ \eta \in \bb X(c) \, | \, \langle \eta | \eta \rangle_{\bb A} = 0 \}$$ and let $\pi^{\bb X}_c$ be the map which quotients by $N^c$.
Note that by Cauchy-Schwarz, $N^c$ is a vector space.
We define $\pi^{\bb X}$ on components by $\pi^{\bb X}(\bb X)(c) = \bb X(c)/N^c$.
Defining it on morphisms takes a bit more work.
Tor $f \colon d \to c$ in $\cC$, there exists a morphism $$\xymatrix{\bb X(c) \ar[r]^{\bb X(f)} & \bb X(d) \ar[r]^<<<<{\pi^{\bb X}_{d}} & \pi^{\bb X}(\bb X)(d)}.$$
By the universal property of the quotient, it suffices to see that $\bb X(f) (N^c) \subset N^d$ to obtain a morphism $\pi^{\bb X}(\bb X)(f) \colon \pi^{\bb X}(\bb X)(c) \to \pi^{\bb X}(\bb X)(d)$.
For $\eta \in \bb X(c)$, we study the diagram

\begin{center}
\begin{tikzpicture}[baseline={(current bounding box.center)}, thick]
\draw (0, 0) node[draw=black,fill=white,very thick,rounded corners,baseline=center] (a){$\langle \cdot | \cdot \rangle_{\bb A}$};
\draw (a.north) --node[pos=1, above] {$\bb A$}++(0, 0.5);
\draw [red, very thick] ([xshift=-0.45cm]a.south) --node[pos=0.4, left, black] {$\bbover X$} ++ (0, -0.75);
\draw [red, very thick] ([xshift=0.45cm]a.south) --node[pos=0.4, right, black] {$\bb X$} ++ (0, -0.75);
\draw ([xshift=-0.45cm, yshift=-1.25cm]a) node[draw=black,fill=white,very thick,rounded corners,baseline=center] (b){$\overline{\eta}$};
\draw ([xshift=0.45cm, yshift=-1.25cm]a) node[draw=black,fill=white,very thick,rounded corners,baseline=center] (c){$\eta$};
\draw (b.south) --node[pos=0.4, left] {$\bbover d$}++ (0, -0.75);
\draw (c.south) --node[pos=0.4, right] {$\bb d$}++ (0, -0.75);
\draw ([yshift=-1.25cm]b) node[draw=black,fill=white,very thick,rounded corners,baseline=center] (d){$\bbover f$};
\draw ([yshift=-1.25cm]c) node[draw=black,fill=white,very thick,rounded corners,baseline=center] (e){$\bb f$};
\draw (d.south) --node[pos=1, below] {$\bbover c$}++ (0, -0.5);
\draw (e.south) --node[pos=1, below] {$\bb c$}++ (0, -0.55);
\end{tikzpicture}
\end{center}

Using associativity and unitality of $\bb A$, and identifying $f$ with the map it is sent to by Frobenius reciprocity, we rewrite this diagram as

\begin{center}
\begin{tikzpicture}[baseline={(current bounding box.center)}, thick]
\draw (0, 0) node[draw=black,fill=white,very thick,rounded corners,baseline=center] (a){$\langle \cdot | \cdot \rangle_{\bb A}$};
\draw (a.north) --node[pos=1, above] {$\bb A$}++(0, 1);
\draw [red, very thick] ([xshift=-0.45cm]a.south) --node[pos=0.4, left, black] {$\bbover X$} ++ (0, -0.75);
\draw [red, very thick] ([xshift=0.45cm]a.south) --node[pos=0.4, right, black] {$\bb X$} ++ (0, -0.75);
\draw ([xshift=-0.45cm, yshift=-1.25cm]a) node[draw=black,fill=white,very thick,rounded corners,baseline=center] (b){$\overline{\eta}$};
\draw (b.south) to [out=-90, in=0] node[pos=0.5, below] {$\bbover d$} ([xshift=-0.6cm, yshift=-0.5cm]b.south) coordinate (f);
\draw ([xshift=0.45cm, yshift=-1.25cm]a) node[draw=black,fill=white,very thick,rounded corners,baseline=center] (c){$\eta$};
\draw (c.south) to [out=-90, in=180] node[pos=0.5, below] {$\bb d$} ([xshift=0.6cm, yshift=-0.5cm]c.south) coordinate (g);
\draw ([xshift=-1.25cm]b) node[draw=black,fill=white,very thick,rounded corners,baseline=center] (d){$\bbover f$};
\draw ([xshift=1.25cm]c) node[draw=black,fill=white,very thick,rounded corners,baseline=center] (e){$\bb f$};
\draw ([xshift=-0.1cm]d.south) --node[pos=1, below] {$\bbover c$}++ (0, -0.75);
\draw ([xshift=0.1cm]e.south) --node[pos=1, below] {$\bb c$}++ (0, -0.8);
\draw ([xshift=0.1cm]d.south) to [out=-90, in=180] node[pos=0.5, below] {$\bb d$} (f);
\draw ([xshift=-0.1cm]e.south) to [out=-90, in=0] node[pos=0.5, below] {$\bbover d$} (g);
\filldraw [black] ([yshift=0.5cm]d.north) circle (2pt) coordinate (h);
\filldraw [black] ([yshift=0.5cm]e.north) circle (2pt) coordinate (i);
\draw (h) to [out=90, in=180] node[pos=0.5, above left] {$\bb A$} ([yshift=0.5cm]a.north);
\draw (i) to [out=90, in=0] node[pos=0.5, above right] {$\bb A$} ([yshift=0.5cm]a.north);
\end{tikzpicture}
\end{center}
Observe that the 4-valent vertex is well-defined by associativity of $\bb A$.
This diagram represents the element $\bbover f \langle \eta | \eta \rangle_{\bb A} \bb f$ of $\bb A(\overline{c} \otimes c)$.
Since $\langle \eta | \eta \rangle_{\bb A} = 0$, $\bbover f \langle \eta | \eta \rangle_{\bb A} \bb f = 0$ as well.
So we obtain a morphism $\pi^{\bb X}(\bb X)(f)$.
By the universal property of the quotient, this morphism is given concretely by $\pi^{\bb X}(\bb X)(f)(\overline{\eta}) = \overline{f(\eta)}$ for $\eta \in \bb X(c)$.
The construction of $\pi^{\bb X}(\bb X)(f)$ makes $\pi^{\bb X}(\bb X)$ into a linear functor.
Now we equip $\pi^{\bb X}(\bb X)$ with a right $\bb A$ action $\pi^{\bb X}(\lambda)$.
To do so, we provide the data of $\pi^{\bb X}(\lambda)_{a, b}$ for $a, b \in \Irr(\cC)$.
Where $\overline{\eta} \in \pi^{\bb X}(\bb X)(a)$ and $\xi \in \bb A(b)$, we define $$\pi^{\bb X}(\lambda)_{a, b}(\overline{\eta} \otimes \xi) = \overline{\lambda_{a, b}(\eta \otimes \xi)}$$ 
We check that this action is well-defined.
For $\overline{\eta} = \overline{\eta'} \in \pi^{\bb X}(\bb X)(a)$ and $\xi \in \bb A(b)$, it suffices to show that $\langle \lambda_{a, b}((\eta - \eta') \otimes \xi) | \lambda_{a, b}((\eta - \eta') \otimes \xi) \rangle_{\bb A} = 0$.
To do so, consider the following diagram, which is manifestly 0 by left and right linearity of $\langle \cdot | \cdot \rangle_{\bb A}$.

\begin{center}
\begin{tikzpicture}[baseline={(current bounding box.center)}, thick]
\draw (0, 0) node[draw=black,fill=white,very thick,rounded corners,baseline=center, inner xsep = 1cm] (a){$\langle \cdot | \cdot \rangle_{\bb A}$};
\draw (a.north) --node[pos=1, above] {$\bb A$}++ (0, 0.25);
\draw [red, very thick] ([xshift=-1cm]a.south) -- node[pos=0.5, right, black] {$\bbover X$} coordinate (b) ++ (0, -1.5);
\draw [red, very thick] ([xshift=1cm]a.south) -- node[pos=0.5, left, black] {$\bb X$} coordinate (c) ++ (0, -1.5);
\draw ([xshift=-1cm, yshift=-2cm]a) node[draw=black,fill=white,very thick,rounded corners,baseline=center] (d){$\overline{\eta - \eta'}$};
\draw ([xshift=1cm, yshift=-2cm]a) node[draw=black,fill=white,very thick,rounded corners,baseline=center] (e){$\eta - \eta'$};
\draw ([xshift=-0.75cm]d.west) node[draw=black,fill=white,very thick,rounded corners,baseline=center] (f){$\overline{\xi}$};
\draw ([xshift=0.75cm]e.east) node[draw=black,fill=white,very thick,rounded corners,baseline=center] (g){$\xi$};
\draw (f.north) to [out = 90, in = 180] node[pos=0.6, above left] {$\bbover A$} ([yshift=0.35cm]b);
\draw (g.north) to [out = 90, in = 0] node[pos=0.6, above right] {$\bb A$} ([yshift=0.35cm]c);
\draw (d.south) --node[pos=1, below] {$\bbover b$} ++ (0, -0.5);
\draw (e.south) --node[pos=1, below] {$\bb b$} ++ (0, -0.5);
\draw (f.south) --node[pos=1, below] {$\bbover a$} ++ (0, -0.5);
\draw (g.south) --node[pos=1, below] {$\bb a$} ++ (0, -0.5);
\end{tikzpicture}
\end{center}

As $\lambda$ is natural, so too is $\pi^{\bb X}(\lambda)$.
To make $\pi^{\bb X}(\bb X)$ into a right $\bb A$ pre-Hilbert $\rmC^*$ module object, we equip $\pi^{\bb X}(\bb X)$ with an $\bb A$-valued inner product.
Where $\widetilde{\eta}$ and $\widetilde{\xi}$ represent equivalence classes in components of $\pi^{\bb X}(\bb X)$, we define $\langle \widetilde{\eta} | \widetilde{\xi} \rangle_{\bb A}^{\pi^{\bb X}(\bb X)} \coloneq \langle \eta | \xi \rangle_{\bb A}^{\bb X}$.
By the Cauchy-Schwarz inequality and since conditional expectations are faithful, the form $\langle \cdot | \cdot \rangle_{\bb A}^{\pi^{\bb X}(\bb X)}$ is a well-defined right $\bb A$-valued inner product on $\pi^{\bb X}(\bb X)$.

Now we complete components with respect to the right $\bb A(1_{\cC})$-valued inner product associated to $\langle \cdot | \cdot \rangle_{\bb A}^{\pi^{\bb X}(\bb X)}$ and call the result $\widetilde{\bb X}$.
We overload notation and denote elements of components of both $\widetilde{\bb X}$ and $\pi^{\bb X}(\bb X)$ without bars.
For $f \colon d \to c $ in $\cC$, we define the morphism $\widetilde{\bb X}(f)$ by $\eta = \lim \eta_n \mapsto \lim \pi^{\bb X}(f)(\eta_n)$.
This is well-defined since limits are unique.
(We omit checking that operations on right module objects whose components are complete are well-defined from here on.)
That $\widetilde{\bb X}$ is a linear functor follows as $\bb X$ is one.
The right $\bb A$ action $\widetilde{\lambda}$ on $\widetilde{\bb X}$ is the right $\bb A$ action $\lambda$ on $\bb X$ extended by continuity.
As $\lambda$ is natural, so too is $\widetilde{\lambda}$.
Thus extending $\langle \cdot | \cdot \rangle_{\bb A}^{\pi^{\bb X}(\bb X)}$ by continuity to a well-defined right $\bb A$-valued inner product $\langle \cdot | \cdot \rangle_{\bb A}^{\widetilde{\bb X}}$ on $\widetilde{\bb X}$ produces a right $\bb A$ Hilbert $\rmC^*$ module object.

By construction, $\pi^{\bb X}$ is a natural transformation $\bb X \to \pi^{\bb X}(\bb X)$.
As completing component-wise is also a natural transformation, constructing $\widetilde{\bb X}$ from $\bb X$ is natural in $\bb X$.
Note that if $\bb X$ is a right $\bb A$ Hilbert $\rmC^*$ module object, then $\widetilde{\bb X} = \bb X$.

Finally, the construction of $\widetilde{\bb X}$ is functorial.
Let $\alpha \colon \bb X \to \bb Y$ be a right $\bb A$-linear morphism, where $\bb X$ and $\bb Y$ are right $\bb A$ module objects with right $\bb A$-valued pre-inner products.
For $c \in \cC$ and $\eta \in \widetilde{\bb X}(c)$, $\eta = \lim \eta_n$ for $\eta_n \in \pi^{\bb X}(\bb X)(c)$.
We let $\widetilde{\alpha}_c(\eta) \coloneq \lim \alpha_c(\eta_n)$.
This defines a right $\bb A$-linear morphism from $\widetilde{\bb X}$ to $\widetilde{\bb Y}$.
Indeed, since $\alpha$ is natural and right $\bb A$-linear, so too is $\widetilde{\alpha}$.
One can quickly check that $\widetilde{1_{\bb X}} = 1_{\widetilde{\bb X}}$ and for $\bb X \stackrel{\alpha}{\to} \bb Y \stackrel{\beta}{\to} \bb Z$ in the category of right $\bb A$ module objects, $\widetilde{\beta \alpha} = \widetilde{\beta} \widetilde{\alpha}$.
\end{const}

Construction \ref{const:completion} satisfies the following universal property.

\begin{prop} 
Let $\bb A$ be a $\rmC^*$ algebra object, $\bb X$ be a right $\bb A$ module object with a right $\bb A$-valued pre-inner product, and let $\bb Y$ be a right $\bb A$ Hilbert $\rmC^*$ module object.
Given a right $\bb A$-linear morphism $\alpha \colon \bb X \to \bb Y$, there exists a unique right $\bb A$-linear morphism $\alpha' \colon \widetilde{\bb X} \to \bb Y$ such that the following diagram commutes
\begin{equation}
\label{eq:UniPropCompletion}
\xymatrix{\bb X \ar[r]^{\widetilde{\cdot}} \ar[dr]_{\alpha} & \widetilde{\bb X} \ar[d]^<<<{\alpha'} \\ & \bb Y}
\end{equation}
\end{prop}

\begin{proof}
For $a \in \cC$ and $\eta = \lim \eta_n \in \widetilde{\bb X}(a)$, let $\alpha'_a(\eta) = \lim \alpha_a \eta_n$.
These assemble into a natural transformation $\alpha'$, as $\alpha$ is natural and $\bb Y(f)$ is continuous.
Moreover, $\alpha'$ is right $\bb A$-linear as $\alpha$ is.
The diagram \ref{eq:UniPropCompletion} commutes by construction.
If $\beta \colon \widetilde{\bb X} \to \bb Y$ also makes the diagram commute, then $\alpha' = \beta$.
Thus the universal property is satisfied by $\alpha'$.
\end{proof}

\end{document}